\documentclass[11pt,french,english]{article}
\usepackage[T1]{fontenc}
\usepackage[latin9]{inputenc}
\usepackage{geometry}
\usepackage{verbatim}
\usepackage{float}
\usepackage{amsmath}
\usepackage{amssymb}
\usepackage{graphicx}
\usepackage{esint}

\makeatletter

\floatstyle{ruled}
\newfloat{algorithm}{tbp}{loa}
\providecommand{\algorithmname}{Algorithm}
\floatname{algorithm}{\protect\algorithmname}

\usepackage{amsfonts}

\usepackage{float}

\usepackage{esint}

\makeatletter

\floatstyle{ruled}
\newfloat{algorithm}{tbp}{loa}
\providecommand{\algorithmname}{Algorithm}
\floatname{algorithm}{\protect\algorithmname}

\usepackage{amsfonts}

\usepackage{subfigure}
\@ifundefined{definecolor}
 {\@ifundefined{definecolor}
 {\usepackage{color}}{}
}{}

\usepackage{algorithmic}
\floatname{algorithm}{Algorithm}

\newtheorem{theorem}{Theorem}

\newtheorem{lemma}[theorem]{Lemma}

\makeatother
\date{}

\makeatother

\usepackage{babel}
\addto\extrasfrench{%
   
   \providecommand{\fg}{\ifdim\lastskip>\z@\unskip\fi~\frqq}
}

\addto\captionsfrench{\renewcommand{\algorithmname}{Algorithme}}

\begin{document}

\title{Distributed Maximum Likelihood for Simultaneous Self-localization
and Tracking in Sensor Networks }

\author{Nikolas~Kantas%
\thanks{N. Kantas is with the Control and Power Group, Department of Electrical
and Electronic Engineering, Imperial College, London, UK, SW7 2AZ,
e-mail: \{n.kantas@imperial.ac.uk\}.%
},~ Sumeetpal~S.~Singh%
\thanks{S.S. Singh is with the Signal Processing lab, Department of Engineering,
University of Cambridge, Trumpington Road, Cambridge, UK, CB2 1PZ,
e-mail: \{sss40@cam.ac.uk\}.%
},~and~Arnaud~Doucet%
\thanks{A. Doucet is with the Department of Statistics, University of Oxford,
1 South Parks Road, Oxford, OX1 3TG, e-mail: doucet@stats.ox.ac.uk %
}}
\maketitle
\begin{abstract}
We show that the sensor self-localization problem can be cast as a
static parameter estimation problem for Hidden Markov Models and we
implement fully decentralized versions of the Recursive Maximum Likelihood
and on-line Expectation-Maximization algorithms to localize the sensor
network simultaneously with target tracking. For linear Gaussian models,
our algorithms can be implemented exactly using a distributed version
of the Kalman filter and a novel message passing algorithm. The latter
allows each node to compute the local derivatives of the likelihood
or the sufficient statistics needed for Expectation-Maximization.
In the non-linear case, a solution based on local linearization in
the spirit of the Extended Kalman Filter is proposed. In numerical
examples we demonstrate that the developed algorithms are able to
learn the localization parameters. 
\end{abstract}
Collaborative tracking, sensor localization, target tracking, maximum
likelihood, sensor networks

\section{Introduction}

This paper is concerned with sensor networks that are deployed to
perform target tracking. A network is comprised of synchronous sensor-trackers
where each node in the network has the processing ability to perform
the computations needed for target tracking. A moving target will
be simultaneously observed by more than one sensor. If the target
is within the field-of-view of a sensor, then that sensor will collect
measurements of the target. Traditionally in tracking a centralized
architecture is used whereby all the sensors transmit their measurements
to a central fusion node, which then combines them and computes the
estimate of the target's trajectory. However, here we are interested
in performing \emph{collaborative tracking,} but without the need
for a central fusion node. Loosely speaking, we are interested in
developing distributed tracking algorithms for networks whose nodes
collaborate by exchanging appropriate messages between neighboring
nodes to achieve the same effect as they would by communicating with
a central fusion node.

A necessary condition for distributed collaborative tracking is that
each node is able to accurately determine the position of its neighboring
nodes in its local frame of reference. (More details in Section \ref{sec:probFormulate}.)
This is essentially an instance of the \emph{self-localization} problem.
In this work we solve the self-localization problem in an on-line
manner. By on-line we mean that self-localization is performed on-the-fly
as the nodes collect measurements of the moving target. In addition,
given the absence of a central fusion node collaborative tracking
and self-localization have to be performed in a fully \emph{decentralized
}manner, which makes necessary the use of message passing between
neighboring nodes.

There is a sizable literature on the self-localization problem. The
topic has been independently pursued by researchers working in different
application areas, most notably wireless communications \cite{Patwari05,plarre2008tracking,Priyantha03,Ihler04,MosesKP03}.
Although all these works tend to be targeted for the application at
hand and differ in implementation specifics, they may however be broadly
summarized into two categories. Firstly, there are works that rely
on direct measurements of distances between neighboring nodes \cite{plarre2008tracking,Priyantha03,Ihler04,MosesKP03}.
The latter is usually estimated from the Received Signal Strength
(RSS) when each node is equipped with a wireless transceiver. Given
such measurements, it is then possible to solve for the geometry of
the sensor network but with ambiguities in translation and rotation
of the entire network remaining. These ambiguities can be removed
if the absolute position of certain nodes, referred to as anchor nodes,
are known. Another approach to self-localization utilizes \emph{beacon}
nodes which have either been manually placed at precise locations,
or their locations are known using a Global Positioning System (GPS).
The un-localized nodes will use the signal broadcast by these beacon
nodes to self-localize \cite{Patwari05,vemula:sensor,cevher06,chen2011sequential}.
We emphasize that in the aforementioned papers self-localization is
performed off-line. The exception is \cite{chen2011sequential}, where
they authors use Maximum Likelihood (ML) and Sequential Monte Carlo
(SMC) in a centralized manner.

In this paper we aim to solve the localization problem without the
need of a GPS or direct measurements of the distance between neighboring
nodes. The method we propose is significantly different.%
{} Initially, the nodes do not know the relative locations of other
nodes, so they can only behave as independent trackers. As the tracking
task is performed on objects that traverse the field of view of the
sensors, information is shared between nodes in a way that allows
them to self-localize. Even though the target's true trajectory is
not known to the sensors, localization can be achieved in this manner
because the same target is being simultaneously measured by the sensors.
This simple fact, which with the exception of \cite{Funiak06,TaylorLaSlat,Baryshnikov07}
seems to have been overlooked in the localization literature, is the
basis of our solution%
\footnote{A short preliminary version of the this work was published in the
conference proceedings \cite{Kantas07}.%
}. However, our work differs from \cite{Funiak06,TaylorLaSlat} in
the application studied as well as the inference scheme. Both \cite{Funiak06,TaylorLaSlat}
formulate the localization as a Bayesian inference problem and approximate
the posterior distributions of interest with Gaussians. \cite{TaylorLaSlat}
uses a moment matching method and appears to be centralized in nature.
The method in \cite{Funiak06} uses instead linearization, is distributed
and on-line, but its implementation relies on communication via a
junction tree (see \cite{Pearl88} for details) and requires an anchor
node as pointed out in \cite[Section 6.2.3]{kantas-thesis}. In this
paper we formulate the sensor localization problem as a static parameter
estimation problem for Hidden Markov Models (HMMs) \cite{cappe2005HMM,elliott1995HMM}
and we estimate these static parameters using a ML approach, which
has not been previously developed for the self-localization problem.
We implement fully \emph{decentralized} versions of the two most common
on-line ML inference techniques, namely Recursive Maximum Likelihood
(RML) \cite{HoL91,LeGlandMevel97,CoR98} and on-line Expectation-Maximization
(EM) \cite{ford-thesis,elliott-online-em,cappe-on-line-em}. A clear
advantage of this approach compared to previous alternatives is that
it makes an on-line implementation feasible. Finally, \cite{Baryshnikov07}
is based on the principle shared by our approach and \cite{Funiak06,TaylorLaSlat}.
In \cite{Baryshnikov07} the authors exploit the correlation of the
measurements made by the various sensors of a hidden spatial process
to perform self-localization. However for reasons concerned with the
applications being addressed, which is not distributed target tracking,
their method is not on-line and is centralized in nature.

In the signal processing literature for sensor networks one may find
various related problems. In \cite{Nowak03} a distributed EM algorithm
was developed to estimate the parameters of a Gaussian mixture used
to model the measurements of a sensor network deployed for environmental
monitoring (see \cite{Sato-ishii} for an on-line version.) In \cite{Blatt04distributedmaximum}
a similar problem is treated using a distributed gradient method.
We emphasize that in each of these papers the measurements correspond
to a static source instead of a dynamically evolving target. In addition,
a related problem is that of \emph{sensor registration,} which aims
to compensate for systematic biases in the sensors and has been studied
by the target tracking community \cite{OkelloChalla03,VermaakMaskellBriers05}.
However, the algorithms devised in \cite{OkelloChalla03,VermaakMaskellBriers05}
are centralized. Yet another related problem is the problem of average
consensus \cite{Xiao05}. The value of a global static parameter is
measured at each node via a linear Gaussian observation model and
the aim is to obtain a maximum likelihood estimate in a distributed
fashion. Note that all the aforementioned papers, except \cite{Funiak06}
and \cite{TaylorLaSlat}, do not deal with a distributed localization
and tracking task.


The structure of the paper is as follows. We begin with the specification
of the statistical model for the localization and tracking problem
in Section \ref{sec:probFormulate}. In Section \ref{sec:distFiltering}
we show how message passing may be utilized to perform distributed
filtering. In Section \ref{sec:distLocalization} we derive the distributed
RML\ and on-line EM algorithms. Section \ref{sec:examples}\ presents
several numerical examples on small and medium sized networks. In
Sections \ref{sec:Discussion} we provide a discussion and a few concluding
remarks. The Appendix contains more detailed derivations of the distributed
versions of RML and EM.

\section{Problem Formulation\label{sec:probFormulate}}

We consider the sensor network $(\mathcal{V},\mathcal{E})$ where
$\mathcal{V}$ denotes the set of nodes of the network and $\mathcal{E}$\ is
the set of edges (or communication links between nodes.) 
We will assume that the sensor network is connected, i.e. for any
pair of nodes $i,j\in\mathcal{V}$ there is at least one path from
$i$ to $j$. Nodes $i,j\in\mathcal{V}$ are adjacent or neighbors
provided the edge $(i,j)\in\mathcal{E}$\ exists. Also, we will assume
that if $(i,j)\in\mathcal{E}$, then $(j,i)\in\mathcal{E}$ as well.
This implies is that communication between nodes is bidirectional.
The nodes observe the same physical target at discrete time intervals
$n\in\mathbb{N}$. We will assume that all sensor-trackers are synchronized
with a common clock and that the edges joining the different nodes
in the network correspond to reliable communication links. These links
define a neighborhood structure for each node and we will also assume
that each sensor can only communicate with its neighboring nodes.

The hidden state, as is standard in target tracking, is defined to
comprise of the position and velocity of the target, $X_{n}^{r}=[X_{n}^{r}(1),X_{n}^{r}(2),X_{n}^{r}(3),X_{n}^{r}(4)]^{\text{T}},$
where $X_{n}^{r}(1)$ and $X_{n}^{r}(3)$\ is the target's $x$ and
$y$ position while $X_{n}^{r}(2)$ and $X_{n}^{r}(4)$\ is the velocity
in the $x$ and $y$ direction. Subscript $n$ denotes time while
superscript $r$ denotes the coordinate system w.r.t. which these
quantities are defined. For generality we assume that each node maintains
a local coordinate system (or frame of reference) and regards itself
as the origin (or center of) its coordinate system.

As a specific example, consider the following linear Gaussian model:
\begin{equation}
X_{n}^{r}=A_{n}X_{n-1}^{r}+b_{n}^{r}+V_{n},\text{\quad}n\geq1,\label{eq:linearGaussianState}
\end{equation}
where\ $V_{n}$\ is zero mean Gaussian additive noise with variance
$Q_{n}$ and $b_{n}^{r}$ are deterministic inputs. The measurement
$Y_{n}^{r}$\ made by node $r$ is also defined relative to the local
coordinate system at node $r$. For a linear Gaussian observation
model the measurement is generated as follows: 
\begin{equation}
Y_{n}^{r}=C_{n}^{r}X_{n}^{r}+d_{n}^{r}+W_{n}^{r},\text{\quad}n\geq1,\label{eq:linearGaussianObs}
\end{equation}
where $W_{n}^{r}$ is zero mean Gaussian additive noise with variance
$R_{n}^{r}$ and $d_{n}^{r}$ is deterministic. Note that the time
varying observation model $\{(C_{n}^{r},d_{n}^{r},R_{n}^{r})\}_{n\geq1}$
is different for each node. A time-varying state and observation model
is retained for an Extended Kalman Filter (EKF) implementation in
the non-linear setting to be defined below. It is in this setting
that the need for sequences $\{b_{n}^{r}\}_{n\geq1}$ and $\{d_{n}^{r}\}_{n\geq1}$
arises. Also, the dimension of the observation vector $Y_{n}^{r}$
need not be the same for different nodes since each node may be equipped
with a different sensor type. For example, node $r$ may obtain measurements
of the target's position while node $v$ measures bearing. Alternatively,
the state-space model in (\ref{eq:linearGaussianState})-(\ref{eq:linearGaussianObs})
can be expressed in the form of a Hidden Markov Model (HMM): 
\begin{align}
X_{n}^{r}|X_{n-1}^{r}=x_{n-1}^{r} & \sim f_{n}(.|x_{n-1}^{r}),\\
Y_{n}^{r}|X_{n}^{r}=x_{n}^{r} & \sim g_{n}^{r}(.|x_{n}^{r}),\label{eq:obsgen}
\end{align}
 where $f_{n}$ denotes the transition density of the target and $g_{n}^{r}$
the density of the likelihood of the observations at each node $r$.

\begin{figure}[h]
\centering

\begin{centering}
\subfigure[Three node joint tracking example]{\includegraphics[width=0.5\textwidth,height=0.5\textwidth]{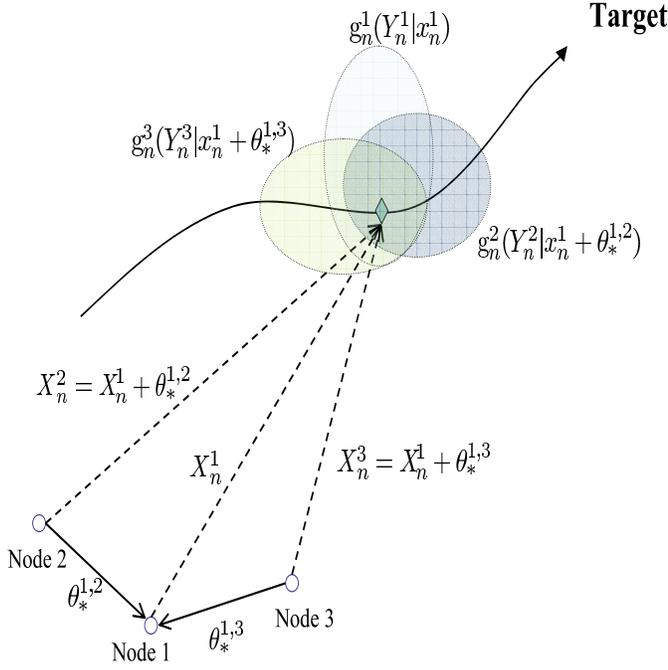}}\subfigure[Joint tracking error vs number of nodes]{
\includegraphics[width=0.5\textwidth,height=0.4\textwidth]{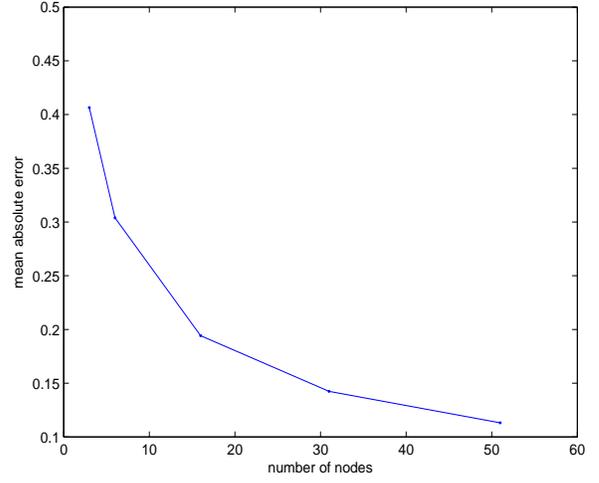}}
\par\end{centering}

\caption{Left: a three node network tracking a target traversing its field
of view. The trajectory of the target is shown with the solid line.
Each node regards itself as the center of its local coordinate system.
At time $n$ a measurement is registered by all three nodes. The ellipses
show the support of the observation densities for the three nodes,
i.e. the support of $g_{n}^{1}(Y_{n}^{1}|.)$ is defined as all $x_{n}^{1}$
such that $g_{n}^{1}(Y_{n}^{1}|x_{n}^{1})>0$ ; similarly for the
rest. The filtering update step at node 1 will clearly benefit from
the observations made by nodes 2 and 3. The localization parameters
$\theta_{*}^{1,2}$, $\theta_{*}^{1,3}$ are the coordinates of node
1 in the local coordinate systems of node 2 and 3 respectively. While
$X_{n}^{r}$ was defined to be the state of the target, which includes
its velocity, for this illustration only, $X_{n}^{r}$ is to be understood
as the position of the target at time $n$ w.r.t. the coordinate system
of node $r$. Right: Average absolute tracking error is plotted against
the number of nodes to illustrate the benefit of collaborative tracking.
The results are obtained using a centralized implementation with 50
independent runs, $10^{4}$ time steps for a chain sensor network
of different length and $A_{n}=B_{n}=Q_{n}=C_{n}^{i}=D_{n}^{i}=R_{n}^{i}=1$,
$b_{n}^{i}=d_{n}^{i}=0$.}

\label{fig:XXX} 
\end{figure}

Figure \ref{fig:XXX} (a) illustrates a three node setting where a
target is being jointly observed and tracked by three sensors. (Only
the position of the target is shown.) At node 1, $X_{n}^{1}$\ is
defined relative to the local coordinate system of node 1 which regards
itself as the origin. Similarly for nodes 2 and 3. We define $\theta_{\ast}^{i,j}$
to be the position of node $i$ \textit{in the local coordinate system}
of node $j$. This means that the vector $X_{n}^{i}$ relates to the
local coordinate system of node $j$ as follows (see Figure \ref{fig:XXX}):
\[
X_{n}^{j}=X_{n}^{i}+\theta_{\ast}^{i,j}.
\]
 The \emph{localization} parameters $\{\theta_{\ast}^{i,j}\}_{(i,j)\in\mathcal{E}}$\ are
static as the nodes are not mobile. We note the following obvious
but important relationship: if nodes $i$ and $j$ are connected through
intermediate nodes $j_{1},j_{2},\ldots,j_{m}$ then 
\begin{equation}
\theta_{\ast}^{i,j}=\theta_{\ast}^{i,j_{1}}+\theta_{\ast}^{j_{1},j_{2}}+\theta_{\ast}^{j_{2},j_{3}}+\ldots+\theta_{\ast}^{j_{m-1},j_{m}}+\theta_{\ast}^{j_{m},j}.\label{eq:sumThetas}
\end{equation}
This relationship is exploited to derive the distributed filtering
and localization algorithms in the next section. We define $\theta_{\ast}^{i,j}$
so that the dimensions are the same as the target state vector. When
the state vector is comprised of the position and velocity of the
target, only the first and third components of $\theta_{\ast}^{i,j}$
are relevant while the other two are redundant and set to $\theta_{\ast}^{i,j}(2)=0$
and $\theta_{\ast}^{i,j}(4)=0$. Let 
\begin{equation}
\theta_{\ast}\equiv\{\theta_{\ast}^{i,j}\}_{(i,j)\in\mathcal{E}},\quad\theta_{\ast}^{i,i}\equiv0,
\end{equation}
 where $\theta_{\ast}^{i,i}$\ for all $i\in\mathcal{V}$ is defined
to be the zero vector.

Let $Y_{n}$ denote all the measurements received by the network at
time $n$, i.e. $Y_{n}\equiv\{Y_{n}^{v}\}_{v\in\mathcal{V}}$. We
also denote the sequence $(Y_{1},...,Y_{n})$ by $Y_{1:n}$. In the
\emph{collaborative or joint} filtering problem, each node $r$ computes
the local filtering density: 
\begin{equation}
p_{\theta_{\ast}}^{r}(x_{n}^{r}|Y_{1:n})\propto p_{\theta_{\ast}}^{r}(Y_{n}|x_{n}^{r})p_{\theta_{\ast}}^{r}(x_{n}^{r}|Y_{1:n-1}),\label{eq:update}
\end{equation}
 where $p_{\theta_{\ast}}^{r}(x_{n}^{r}|Y_{1:n-1})$ is the predicted
density and is related to the filtering density of the previous time
through the following prediction step: 
\begin{equation}
p_{\theta_{\ast}}^{r}(x_{n}^{r}|Y_{1:n-1})=\int f_{n}(x_{n}^{r}|x_{n-1}^{r})p_{\theta_{\ast}}^{r}(x_{n-1}^{r}|Y_{1:n-1})dx_{n-1}^{r}.\label{eq:prediction}
\end{equation}
 The likelihood term is 
\begin{equation}
p_{\theta_{\ast}}^{r}(Y_{n}|x_{n}^{r})=\prod\limits _{v\in\mathcal{V}}g_{n}^{v}(Y_{n}^{v}|x_{n}^{r}+\theta_{\ast}^{r,v}),\label{eq:collabLikelihood}
\end{equation}
where the superscript on the densities indicate the coordinate system
they are defined w.r.t. (and the node the density belongs to) while
the subscript makes explicit the dependence on the localization parameters.
Let also $\mu_{n|n-1}^{r}$ and $\mu_{n}^{r}$ denote the predicted
and filtered mean of the densities $p_{\theta_{*}}^{r}(x_{n}^{r}|Y_{1:n-1})$
and $p_{\theta_{*}}^{r}(x_{n}^{r}|Y_{1:n})$ respectively, where the
dependence on $\theta_{*}$ is suppressed in the notation. The prediction
step in (\ref{eq:prediction}) can be implemented locally at each
node without exchange of information, but the update step in (\ref{eq:update})
incorporates all the measurements of the network. Figure \ref{fig:XXX}
(a) shows the support of the three observation densities as ellipses
where the support of $g_{n}^{1}(Y_{n}^{1}|\cdot)$ is defined to be
all $x^{1}$\ such that $g_{n}^{1}(Y_{n}^{1}|\cdot)>0$; similarly
for the rest. The filtering update step at node 1 can only include
the observations made by nodes 2 and 3 provided the localization parameters
$\theta_{\ast}^{1,2}$\ and $\theta_{\ast}^{1,3}$\ are known locally
to node 1, since the likelihood $p_{\theta_{\ast}}^{1}(Y_{n}|x_{n}^{1})$
defined in (\ref{eq:collabLikelihood}) is 
\[
g_{n}^{1}(Y_{n}^{1}|x_{n}^{1})g_{n}^{2}(Y_{n}^{2}|x_{n}^{1}+\theta_{\ast}^{1,2})g_{n}^{3}(Y_{n}^{3}|x_{n}^{1}+\theta_{\ast}^{1,3}).
\]

The term joint filtering is used since each sensor benefits from the
observation made by all the other sensors. An illustration of the
benefit w.r.t. the tracking error is in Figure \ref{fig:XXX} (b).
We will show in Section \ref{sec:distFiltering} that it is possible
to implement joint filtering in a truly distributed manner, i.e. each
node executes a message passing algorithm (with communication limited
only to neighboring nodes) that is scalable with the size of the network.
However joint filtering hinges on knowledge of the localization parameters
$\theta_{\ast}$ which are unknown \emph{a priori}. In Section \ref{sec:distLocalization}
we will propose distributed estimation algorithms to learn the localization
parameters, which refine the parameter estimates as new data arrive.
These proposed algorithms in this context are to the best of our knowledge
novel.

\subsection{Non-linear Model \label{sub:Non-linear-Model}}

Most tracking problems of practical interest are essentially non-linear
non-Gaussian filtering problems. SMC methods, also known as Particle
Filters, provide very good approximations to the filtering densities
\cite{DFG01Book}. While it is possible to develop SMC methods for
the problem presented here, the resulting algorithms require significantly
higher computational cost. We refer the interested reader to \cite[Chapter 9]{kantas-thesis}
for more details. In the interest of execution speed and simplicity,
we employ the linearization procedure of the Extended Kalman filter
(EKF) when dealing with a non-linear system. Specifically, let the
distributed tracking system be given by the following model: 
\begin{align}
X_{n}^{r} & =\phi_{n}(X_{n-1}^{r})+V_{n},\\
Y_{n}^{r} & =\psi_{n}^{r}(X_{n}^{r})+W_{n}^{r},
\end{align}
where $\phi_{n}:\mathbb{R}^{4}\rightarrow\mathbb{R}^{4}$ and $\psi_{n}^{r}:\mathbb{R}^{4}\rightarrow\mathbb{R}^{d_{y}}$
are smooth continuous functions. At time $n$, 
each node will linearize its state and observation model about the
filtered and predicted mean respectively. Specifically, a given node
$r$ will implement: 
\begin{align}
X_{n}^{r}=\phi_{n}(\mu_{n-1}^{r}) & +\nabla\phi_{n}(\mu_{n-1}^{r})(X_{n-1}^{r}-\mu_{n-1}^{r})+V_{n},\\
Y_{n}^{r}=\psi_{n}^{r}(\mu_{n|n-1}^{r}) & +\nabla\psi_{n}^{r}(\mu_{n|n-1}^{r})(X_{n}^{r}-\mu_{n|n-1}^{r})+W_{n}^{r}.
\end{align}
where for a mapping $f:\mathbb{R}^{d}\rightarrow\mathbb{R}^{d}$,
$\nabla f\equiv\lbrack\nabla f_{1},\ldots,\nabla f_{d}]^{\text{T}}$.
Note that after linearization extra additive terms appear as seen
in the setting described by equations (\ref{eq:linearGaussianState})-(\ref{eq:linearGaussianObs}).

\subsection{Message passing}

Assume at time $n$, the estimate of the localization parameters is
$\theta_{n}=\{\theta_{n}^{i,j}\}_{(i,j)\in\mathcal{E}}$, with $\theta_{n}^{i,j}$
known to node $j$ only. To perform the prediction and update steps
in (\ref{eq:update})-(\ref{eq:prediction}) locally at each node
a naive approach might require each node to access to all localization
parameters $\theta_{n}$ and all the different model parameters $\{(C_{n}^{r},d_{n}^{r},R_{n}^{r})\}_{n\geq1,r\in\mathcal{V}}$
. A scheme that requires all this information to be passed at every
node would be inefficient. It would require a prohibitive amount of
communication even for relatively few nodes and redundant computations
would be performed at the different nodes. The core idea in this paper
is to avoid this by storing the parameters in $\theta_{n}$ across
the network and perform required computations only at the nodes where
the parameters are stored. The results of these computations are then
propagated in the network using an efficient message passing scheme.

\begin{algorithm}
\selectlanguage{french}%
\begin{algorithmic}[1]

\STATE \textbf{begin}\foreignlanguage{english}{ }

\STATE \foreignlanguage{english}{At $k=1,$ compute:
\begin{align}
m_{n,1}^{i,j} & =F_{n}^{i},\\
\ddot{m}_{n,1}^{i,j} & =F_{n}^{i}\theta_{n}^{j,i}.
\end{align}
}\STATE  \textbf{for}\foreignlanguage{english}{ $k=2,...,K$ compute:}

\selectlanguage{english}%
\begin{align}
m_{n,k}^{i,j} & =F_{n}^{i}+\sum\limits _{p\in\text{ne}(i)\setminus\{j\}}m_{n,k-1}^{p,i},\\
\ddot{m}_{n,k}^{i,j} & =m_{n,k}^{i,j}\theta_{n}^{j,i}+\sum\limits _{p\in\text{ne}(i)\setminus\{j\}}\ddot{m}_{n,k-1}^{p,i}.
\end{align}
\foreignlanguage{french}{\STATE \textbf{endfor}}

\selectlanguage{french}%
\STATE \textbf{end}

\end{algorithmic}\foreignlanguage{english}{\caption{Generic message passing at time $n$}
\label{alg:Messages}}\selectlanguage{english}
\end{algorithm}

Message passing is an iterative procedure with $k=1,\ldots,K$ iterations
for each time $n$ and is steered towards the development of a distributed
Kalman filter, whose presentation is postponed for the next section.
In Algorithm \ref{alg:Messages} we define a recursion of messages
which are to be communicated between all pairs of neighboring nodes
in both directions. Here ne$(i)$ denote the neighbors of node $i$
excluding node $i$ itself. At iteration $k$ the computed messages
from node $i$ to $j$ are matrix and vector quantities of appropriate
dimensions and are denoted as $m_{n,k}^{i,j}$ and $\ddot{m}_{n,k}^{i,j}$
respectively. The source node is indicated by the first letter of
the superscript. Note that during the execution of Algorithm \ref{alg:Messages}
time $n$ remains fixed and iteration $k$ should not be confused
with time $n$. Clearly we assume that the sensors have the ability
to communicate much faster than collecting measurements. 
We proceed with a simple (but key) lemma concerning the aggregations
of sufficient statistics locally at each node.

\begin{lemma} \label{lem:messages} At time $n$, let $\left\{ F_{n}^{v}\right\} _{v\in\mathcal{V}}$
be a collection of matrices where $F_{n}^{v}$ is known to node $v$
only, and consider the task of computing $\sum_{v\in\mathcal{V}}F_{n}^{v}$\ and
$\sum_{v\in\mathcal{V}}F_{n}^{v}\theta_{n}^{r,v}$ at each node $r$
of a network with a tree topology. Using Algorithm 1 and if $K$ is
at least as large as the number of edges connecting the two farthest
nodes in the network, then $\sum_{v\in\mathcal{V}}F_{n}^{v}=F_{n}^{r}+\sum\limits _{j\in\text{ne}(r)}m_{n,K}^{j,r}$
and $\sum_{v\in\mathcal{V}}F_{n}^{v}\theta_{n}^{r,v}=\sum\limits _{j\in\text{ne}(r)}\ddot{m}_{n,K}^{j,r}$.
\end{lemma}

(The proof, which uses (\ref{eq:sumThetas}), is omitted.) %
An additional advantage here is that if the network is very large,
in the interest of speed one might be interested in settling with
computing the presented sums only for a subset of nodes and thus use
a smaller $K$. This also applies when a target traverses the field
of view of the sensors swiftly and is visible only by few nodes at
each time. Finally, a lower value for $K$ is also useful when cycles
are present in order to avoid summing each $F_{n}^{i}$ more than
once, albeit summing only over a subset of $\mathcal{V}$.

\section{Distributed Joint Filtering \label{sec:distFiltering}}

For a linear Gaussian system, the joint filter $p_{\theta}^{r}(x_{n}^{r}|Y_{1:n})$
at node $r$ is a Gaussian distribution with a specific mean vector
$\mu_{n}^{r}$\ and covariance matrix $\Sigma_{n}^{r}$. The derivation
of the Kalman filter to implement $p_{\theta}^{r}(x_{n}^{r}|Y_{1:n})$
is standard upon noting that the measurement model at node $r$ can
be written as $Y_{n}=C_{n}X_{n}^{r}+d_{n}+W_{n}$ where the $i$-th
block of $Y_{n}$, $Y_{n}^{i}$, satisfies $Y_{n}^{i}=C_{n}^{i}(X_{n}^{r}+\theta^{r,i})+d_{n}^{i}+W_{n}^{i}$.
However, there will be {}``non-local'' steps due to the requirement
that quantities $\underset{i\in\mathcal{V}}{\sum}(C_{n}^{i})^{\text{T}}(R_{n}^{i})^{-1}C_{n}^{i}$,
$\underset{i\in\mathcal{V}}{\sum}(C_{n}^{i})^{\text{T}}(R_{n}^{i})^{-1}Y_{n}^{i}$
and $\underset{i\in\mathcal{V}}{\sum}(C_{n}^{i})^{\text{T}}(R_{n}^{i})^{-1}C_{n}^{i}\theta^{r,i}$
be available locally at node $r$. To solve this problem, we may use
Lemma \ref{lem:messages} with $F_{n}^{i}=(C_{n}^{i})^{\text{T}}(R_{n}^{i})^{-1}C_{n}^{i}$
and in order to compute $\underset{i\in\mathcal{V}}{\sum}(C_{n}^{i})^{\text{T}}(R_{n}^{i})^{-1}Y_{n}^{i}$
we will define $\dot{m}_{n,k}^{i,j}$ that is an additional message
similar to $m_{n,k}^{i,j}$.

Recall that $b_{n}^{i},d_{n}^{i}$ are known local variables that
arose due to linearization. Also to aid the development of the distributed
on-line localization algorithms in Section \ref{sec:distLocalization},
we assume that for the time being the localization parameter estimates
$\{\theta_{n}\}_{n\geq1}$ are time-varying and known to the relevant
nodes they belong. For the case where that $b_{n}^{i},d_{n}^{i}=0$,
we summarize the resulting distributed Kalman filter in Algorithm
\ref{alg:distFilter}, which is to be implemented at every node of
the network. Note that messages (\ref{eq:message1})-(\ref{eq:message3})
are matrix and vector valued quantities and require a fixed amount
of memory regardless of the number of nodes in the network. Also,
the same rule for generating and combining messages are implemented
at each node. The distributed Kalman filter presented here bears a
similar structure to the one found in \cite{durrant-whyte94}. However,
the message passing scheme is different and due to the application
in mind we have extra terms relevant to the localization parameters.

\begin{algorithm}
\selectlanguage{french}%
\begin{algorithmic}[1]

\STATE  \textbf{begin }

\STATE \textbf{for}\foreignlanguage{english}{ $n\geq1$: }

\STATE \foreignlanguage{english}{Let the localization parameter be
$\theta_{n}$ and the set of collected measurements be $Y_{n}=\{Y_{n}^{v}\}_{v\in\mathcal{V}}$.
Initialize messages $(m_{n,k}^{i,j},\dot{m}_{n,k}^{i,j},\ddot{m}_{n,k}^{i,j})$
and $(m_{n,k}^{j,i},\dot{m}_{n,k}^{j,i},\ddot{m}_{n,k}^{j,i})$ for
all neighboring nodes $(i,j)\in\mathcal{E}$ as:
\begin{align*}
m_{n,1}^{i,j} & =(C_{n}^{i})^{\text{T}}(R_{n}^{i})^{-1}C_{n}^{i},\\
\dot{m}_{n,1}^{i,j} & =(C_{n}^{i})^{\text{T}}(R_{n}^{i})^{-1}Y_{n}^{i},\\
\ddot{m}_{n,1}^{i,j} & =m_{n}^{i,j}\theta_{n}^{j,i},
\end{align*}
}

\STATE  \textbf{for $k=2,\ldots,K$}\foreignlanguage{english}{ exchange
the messages $(m_{n,k}^{i,j},\dot{m}_{n,k}^{i,j},\ddot{m}_{n,k}^{i,j})$
and $(m_{n,k}^{j,i},\dot{m}_{n,k}^{j,i},\ddot{m}_{n,k}^{j,i})$ defined
below between all neighboring nodes $(i,j)\in\mathcal{E}$:
\begin{align}
m_{n,k}^{i,j} & =(C_{n}^{i})^{\text{T}}(R_{n}^{i})^{-1}C_{n}^{i}+\sum\limits _{p\in\text{ne}(i)\setminus\{j\}}m_{n,k-1}^{p,i},\label{eq:message1}\\
\dot{m}_{n,k}^{i,j} & =(C_{n}^{i})^{\text{T}}(R_{n}^{i})^{-1}Y_{n}^{i}+\sum\limits _{p\in\text{ne}(i)\setminus\{j\}}\dot{m}_{n,k-1}^{p,i},\label{eq:message2}\\
\ddot{m}_{n,k}^{i,j} & =m_{n}^{i,j}\theta_{n}^{j,i}+\sum\limits _{p\in\text{ne}(i)\setminus\{j\}}\ddot{m}_{n,k-1}^{p,i},\label{eq:message3}
\end{align}
}\STATE  \textbf{end for}

\STATE \foreignlanguage{english}{Update the local filtering densities
at each node $r\in\mathcal{V}$: 
\begin{align}
\mu_{n|n-1}^{r} & =A_{n}\mu_{n-1}^{r},\quad\Sigma_{n|n-1}^{r}=A_{n}\Sigma_{n-1}^{r}A_{n}^{\text{T}}+Q_{n},\label{eq:distrKalmanSecond}\\
M_{n}^{r} & =(\Sigma_{n|n-1}^{r})^{-1}+(C_{n}^{r})^{\text{T}}(R_{n}^{r})^{-1}C_{n}^{r}+\sum\limits _{i\in\text{ne}(r)}m_{n}^{i,r}\\
z_{n}^{r} & =(\Sigma_{n|n-1}^{r})^{-1}\mu_{n|n-1}^{r}+(C_{n}^{r})^{\text{T}}(R_{n}^{r})^{-1}Y_{n}^{r}\label{eq:distrKalmanSecondLast}\\
 & +\sum\limits _{i\in\text{ne}(r)}\left(\dot{m}_{n}^{i,r}-\ddot{m}_{n}^{i,r}\right),\notag\nonumber \\
\Sigma_{n}^{r} & =(M_{n}^{r})^{-1},\quad\mu_{n}^{r}=\Sigma_{n}^{r}z_{n}^{r},\label{eq:distrKalmanLast}
\end{align}
}

\STATE  \textbf{end for}

\STATE  \textbf{end}

\end{algorithmic}

\selectlanguage{english}%
\caption{\textbf{Distributed Filtering}}
\label{alg:distFilter}
\end{algorithm}

In the case $b_{n}^{i},d_{n}^{i}\neq0$ modifications to Algorithm
\ref{alg:distFilter} are as follows: in (\ref{eq:distrKalmanSecond}),
to the right hand side of $\mu_{n|n-1}^{r}$, the term $b_{n}^{r}$
should be added and all instances of $Y_{n}^{r}$ should be replaced
with $Y_{n}^{r}-d_{n}^{r}$. Therefore the assuming $b_{n}^{i},d_{n}^{i}=0$
does not compromise the generality of the approach. A direct application
of this modification is the distributed EKF, which is obtained by
adding the term $\phi_{n}(\mu_{n-1}^{r})-\nabla\phi_{n}(\mu_{n-1}^{r})\mu_{n-1}^{r}$
to the right hand side of $\mu_{n|n-1}^{r}$ in (\ref{eq:distrKalmanSecond}),
and replacing all instances of $Y_{n}^{r}$ with $Y_{n}^{r}-\psi_{n}^{r}(\mu_{n|n-1}^{r})+\nabla\psi_{n}^{r}(\mu_{n|n-1}^{r})\mu_{n|n-1}^{r}$.
In addition, one needs to replace $A_{n}$ with $\nabla\phi_{n}(\mu_{n-1}^{r})$.

\section{Distributed Collaborative Localization\label{sec:distLocalization}}

Following the discussion in Section \ref{sec:probFormulate} we will
treat the sensor localization problem as a static parameter estimation
problem for HMMs. The purpose of this section is to develop a fully
decentralized implementation of popular Maximum Likelihood (ML) techniques
for parameter estimation in HMMs. We will focus on two on-line ML
estimation methods: Recursive Maximum Likelihood (RML) and Expectation-Maximization
(EM). For the sake of completeness, we have added brief descriptions
of these techniques in Section \ref{sub:MLE} of the appendix. 


The core idea in our distributed ML formulation is to store the parameter
$\theta_{n}=\{\theta_{n}^{i,j}\}_{(i,j)\in\mathcal{E}}$ across the
network. Each node $r$ will use the available data $Y_{1:n}$ from
every node to estimate $\theta_{\ast}^{r,j}$, which is the component
of $\theta_{*}$ corresponding to edge $(r,j)$. This can be achieved
computing at each node $r$ the ML estimate: 
\begin{equation}
\widetilde{\theta}_{n}^{r,j}=\arg\max_{\theta^{r,j}\in\mathbb{R}^{4}}\log p_{\theta}^{r}(Y_{1:n}).\label{eq:distributedMLestim}
\end{equation}
 Note that each node maximizes its {}``local'' likelihood function
although all the data across the network is being used.

On-line parameter estimation techniques like the RML and on-line EM
are suitable for sensor localization in surveillance applications
because we expect a practically indefinite length of observations
to arrive sequentially. For example, objects will persistently traverse
the field of view of these sensors, i.e. the departure of old objects
would be replenished by the arrival of new ones. A recursive procedure
is essential to give a quick up-to-date parameter estimate every time
a new set of observations is collected by the network. This is done
by allowing every node $r$ to update the estimate of the parameter
along edge $(r,j)$, $\theta_{n}^{r,j}$, according to a rule like
\begin{equation}
\theta_{n+1}^{r,j}=G_{n+1}^{r,j}(\theta_{n},Y_{n}),\qquad n\geq1,\label{eq:titsRML-1}
\end{equation}
where $G_{n+1}^{r,j}$ is an appropriate function to be defined. Similarly
each neighbor $j$ of $r$ will perform a similar update along the
same edge only this time it will update $\theta_{n}^{j,r}$. While
updating both parameters associated to each edge is redundant, it
allows a fully decentralized implementation since no other communication
is needed other than the messages defined in Algorithm \ref{alg:Messages}.
Alternatively one could assign both parameters of an edge to just
one controlling node. For example in the three node network of Figure
\ref{fig:XXX}, the parameters of edge $(1,2)$, $\theta_{n}^{1,2}$
and $\theta_{n}^{2,1}$, could be assigned to node $2$, with the
latter having at each time $n$ to update $\theta_{n}^{2,1}$ using
an expression like (\ref{eq:titsRML-1}) and then send $\theta_{n}^{1,2}=-\theta_{n}^{2,1}$
to node $1$. %

\subsection{Distributed RML\label{sec:distRML}}

For distributed RML, each node $r$ updates the parameter of edge
$(r,j)$ using 
\begin{equation}
\theta_{n+1}^{r,j}=\theta_{n}^{r,j}+\gamma_{n+1}^{r}\left[\nabla_{\theta^{r,j}}\log\int p_{\theta}^{r}(Y_{n}|x_{n}^{r})p_{\theta}^{r}(x_{n}^{r}|Y_{1:n-1})dx_{n}^{r}\right]_{\theta=\theta_{n}},\label{eq:almostOnlineDistributedRML}
\end{equation}
 where $\gamma_{n+1}^{r}$ is a step-size that should satisfy $\sum_{n}{\gamma_{n}^{r}}=\infty$
and $\sum_{n}{\left(\gamma_{n}^{r}\right)}^{2}<\infty$.

The gradient in (\ref{eq:almostOnlineDistributedRML}) is w.r.t. $\theta^{r,j}$.
The local joint \textit{predicted} density $p_{\theta}^{r}(x_{n}^{r}|Y_{1:n-1})$
at node $r$ was defined in (\ref{eq:prediction}) and is a function
of $\theta=\{\theta^{i,j}\}_{(i,j)\in\mathcal{E}}$, and likelihood
term is given in (\ref{eq:collabLikelihood}). Also, the gradient
is evaluated at $\theta_{n}=\{\theta_{n}^{i,j}\}_{(i,j)\in\mathcal{E}}$
while only $\theta_{n}^{r,j}$\ is available locally at node $r$.
The remaining values $\theta_{n}$ are stored across the network.
All nodes of the network will implement such a local gradient algorithm
with respect to the parameter associated to its adjacent edge. We
note that (\ref{eq:almostOnlineDistributedRML}) in the present form
is not an on-line parameter update like (\ref{eq:titsRML-1}) as it
requires browsing through the entire history of observations. This
limitation is removed by defining certain intermediate quantities
that facilitate the online evaluation of this gradient in the spirit
of \cite{LeGlandMevel97,CoR98} (see in the Appendix for more details).

\begin{algorithm}
\selectlanguage{french}%
\begin{algorithmic}[1]

\selectlanguage{english}%
\STATE\foreignlanguage{french}{\textbf{begin}} 

\STATE\foreignlanguage{french}{\textbf{for}} $n\geq1$: let the current
parameter estimate be $\theta_{n}$. Upon obtaining measurements $Y_{n}=\{Y_{n}^{v}\}_{v\in\mathcal{V}}$
the following filtering and parameter update steps are to be performed.

\STATE\textbf{Filtering step}: Perform steps (3-6) in Algorithm \ref{alg:distFilter}.

\STATE\textbf{Parameter update}: Each node $r\in\mathcal{V}$ of
the network will update the following quantities for every edge $(r,j)\in\mathcal{E}$:
\begin{align}
\dot{\mu}_{n|n-1}^{r,j} & =A_{n}\dot{\mu}_{n-1}^{r,j},\label{eq:rml1}\\
\dot{z}_{n}^{r,j} & =(\Sigma_{n|n-1}^{r})^{-1}\dot{\mu}_{n|n-1}^{r,j}-m_{n,K}^{j,r},\\
\dot{\mu}_{n}^{r,j} & =(M_{n}^{r})^{-1}\dot{z}_{n}^{r,j}.\label{eq:rml3}
\end{align}
 Upon doing so the localization parameter is updated: 
\begin{align*}
\theta_{n+1}^{r,j} & =\theta_{n}^{r,j}+\gamma_{n+1}^{r}[-(\dot{\mu}_{n|n-1}^{r,j})^{\text{T}}(\Sigma_{n|n-1}^{r})^{-1}\mu_{n|n-1}^{r}\\
 & +(\dot{z}_{n}^{r,j})^{\text{T}}(M_{n}^{r})^{-1}z_{n}^{r}+\dot{m}_{n,K}^{j,r}-\ddot{m}_{n,K}^{j,r}].
\end{align*}
\STATE\foreignlanguage{french}{\textbf{end for}}

\STATE\foreignlanguage{french}{\textbf{end}}

\selectlanguage{french}%
\end{algorithmic}

\selectlanguage{english}%
\caption{\textbf{Distributed RML}}
\label{alg:rml} 
\end{algorithm}

The distributed RML implementation for self-localization and tracking
is presented in Algorithm \ref{alg:rml}, while the derivation of
the algorithm is presented in the Appendix. The intermediate quantities
(\ref{eq:rml1})-(\ref{eq:rml3}) take values in $\mathbb{R}^{4\times2}$
and may be initialized to zero matrices. For the non-linear model,
when an EKF implementation is used for Algorithm \ref{alg:distFilter},
then Algorithm \ref{alg:rml} remains the same.

\subsection{Distributed on-line EM}

\label{sec:distEM}

We begin with a brief description of distributed EM in an off-line
context and then present its on-line implementation. Given a batch
of $T$\ observations, let $p$ be the (off-line) iteration index
and $\theta_{p}=\{\theta_{p}^{i,j}\}_{(i,j)\in\mathcal{E}}$ be the
current estimate of $\theta_{\ast}$ after $p-1$ distributed EM iterations
on the batch of observations $Y_{1:T}$. Each edge controlling node
$r$ will execute the following E and M steps to update the estimate
of the localization parameter for its edge. For iteration $p=1,2,\ldots$
\begin{align*}
\mbox{(E step) }Q^{r}(\theta_{p},\theta) & =\int\log p_{\theta}^{r}(x_{1:T}^{r},Y_{1:T})p_{\theta_{p}}^{r}(x_{1:T}^{r}|Y_{1:T})dx_{1:T}^{r},\\
\mbox{(M step)}\qquad\theta_{p+1}^{r,j} & =\arg\underset{\theta^{r,j}}{\max}\quad Q^{r}(\theta_{p},(\theta^{r,j},\theta_{p}^{-(r,j)})),
\end{align*}
 where $\theta_{p}^{-(r,j)}=\{\theta_{p}^{e}\}_{e\in\mathcal{E}\backslash(r,j)}$.

To show how the E-step can be computed we write $p_{\theta}^{r}(x_{1:T}^{r},Y_{1:T})$\ as,
\[
p_{\theta}^{r}(x_{1:T}^{r})p_{\theta}^{r}(Y_{1:T}|x_{1:T}^{r})=\prod\limits _{n=1}^{T}f_{n}(x_{n}^{r}|x_{n-1}^{r})p_{\theta}^{r}(Y_{n}|x_{n}^{r}),
\]
 where $p_{\theta}^{r}(Y_{n}|x_{n}^{r})$\ was defined in (\ref{eq:collabLikelihood}).
Note that $p_{\theta_{p}}^{r}(x_{1:T}^{r}|Y_{1:T})$ is a function
of $\theta_{p}=\{\theta_{p}^{i,i^{\prime}}\}_{(i,i^{\prime})\in\mathcal{E}}$\ (and
not just $\theta_{p}^{r,j}$) and the $\theta$-dependance of $p_{\theta}^{r}(x_{1:T}^{r},Y_{1:T})$\ arises
through the likelihood term only as $p_{\theta}^{r}(x_{1:T}^{r})$\ is
$\theta$-independent. This means that in order to compute the E-step,
it is sufficient to maintain the smoothed marginals: 
\[
p_{\theta}^{r}(x_{n}^{r}|Y_{1:T})\propto\int p_{\theta}^{r}(x_{1:T}^{r},Y_{1:T})dx_{1:T\backslash\{n\}}^{r},
\]
where $1\leq n\leq T$ and $dx_{1:T\backslash\{n\}}^{r}$\ means
integration w.r.t. all variables except $x_{n}^{r}$. For linear Gaussian
models this smoothed density is also Gaussian, with its mean and covariance
denoted by $\mu_{n|T}^{r},\Sigma_{n|T}^{r}$ respectively.

The M-step is solved by setting the derivative of $Q^{r}(\theta_{p},(\theta^{r,j},\theta_{p}^{-(r,j)}))$\ w.r.t.
$\theta^{r,j}$ to zero. The details are presented in the Appendix
and the main result is: 
\[
\nabla_{\theta^{r,j}}\int\log p_{\theta}^{r}(Y_{n}|x_{n}^{r})p_{\theta_{p}}^{r}(x_{n}^{r}|Y_{1:T})dx_{n}^{r}=\dot{m}_{n,K}^{j,r}-\ddot{m}_{n,K}^{j,r}-(m_{n,K}^{j,r})^{\text{T}}\mu_{n|T}^{r},
\]
 where $(m_{n,K}^{j,r},\dot{m}_{n,K}^{j,r},\ddot{m}_{n,K}^{j,r}),$
defined in (\ref{eq:message1})-(\ref{eq:message3}), are propagated
with localization parameter $\theta_{p}$ for all observations from
time $1$ to $T$. Only $\ddot{m}_{n,K}^{j,r}$\ is a function of
$\theta^{r,j}$. To perform the M-step, the following equation is
solved for $\theta^{r,j}$ 
\begin{equation}
(\sum_{n=1}^{T}m_{n,K}^{j,r})\theta^{r,j}=\sum_{n=1}^{T}(\dot{m}_{n,K}^{j,r}-(m_{n,K}^{j,r})^{\text{T}}\mu_{n|T}^{r}-\ddot{m}_{n,K}^{j,r}+\ddot{m}_{n,1}^{j,r}).\label{eq:mstep}
\end{equation}
Note that $\theta^{r,j}$ is a function of quantities available locally
to node $r$ only. The M-step can also be written as the following
function: 
\[
\Lambda(\mathcal{S}_{T,1}^{r,j},\mathcal{S}_{T,2}^{r,j},\mathcal{S}_{T,3}^{r,j})=\left(\mathcal{S}_{T,2}^{r,j}\right)^{-1}\left(\mathcal{S}_{T,3}^{r,j}-\mathcal{S}_{T,1}^{r,j}\right),
\]
 where $\mathcal{S}_{T,1}^{r,j}$, $\mathcal{S}_{T,2}^{r,j}$, $\mathcal{S}_{T,3}^{r,j}$
are three summary statistics of the form:
\[
\mathcal{S}_{T,l}^{r,j}=\frac{1}{T}\int\left(\sum_{n=1}^{T}s_{n,l}^{r,j}(x_{n}^{r},Y_{n})\right)p_{\theta_{p}}^{r}(x_{n}^{r}|Y_{1:T})dx_{n}^{r},\: l=1,2,3,
\]
with $s_{n,l}^{r,j}$ being defined as follows: 
\begin{eqnarray*}
s_{n,1}^{r,j}(x_{n}^{r},Y_{n}) & = & (m_{n,K}^{j,r})^{\text{T}}x_{n}^{r},\qquad s_{n,2}^{r,j}(x_{n}^{r},Y_{n})=m_{n,K}^{j,r},\\
\: s_{n,3}^{r,j}(x_{n}^{r},Y_{n}) & = & \dot{m}_{n,K}^{j,r}-\ddot{m}_{n,K}^{j,r}+\ddot{m}_{n,1}^{j,r}.
\end{eqnarray*}
Note that for this problem $s_{n,2}^{r,j}$ and $s_{n,3}^{r,j}$ are
state independent.

An on-line implementation of EM follows by computing recursively running
averages for each of the three summary statistics, which we will denote
as $\mathcal{S}_{n,1}^{r,j},\mathcal{S}_{n,2}^{r,j},\mathcal{S}_{n,3}^{r,j}$.
At each time $n$ these will be used at every node $r$ to update
$\theta^{r,j}$ using $\theta_{n+1}^{r,j}=\Lambda(\mathcal{S}_{n,1}^{r,j},\mathcal{S}_{n,2}^{r,j},\mathcal{S}_{n,3}^{r,j})$.
Note that $\Lambda$ is the same function for every node. The on-line
implementation of distributed EM is found in Algorithm \ref{alg:onlineEM}.
All the steps are performed with quantities available locally at node
$r$ using the exchange of messages as detailed in Algorithm \ref{alg:distFilter}.
The derivation of the recursions for $\mathcal{S}_{n,1}^{r,j}$, $\mathcal{S}_{n,2}^{r,j}$,
$\mathcal{S}_{n,3}^{r,j}$ are based on (\ref{eq:updateStatonline})-(\ref{eq:updateStatonline2})
in the Appendix. Here $\gamma_{n}^{r}$ is a step-size satisfying
the same conditions as in RML and $\theta_{0}$ can be initialized
arbitrarily, e.g. the zero vector.  Finally, it has been reported
in \cite{cappe-on-line-em-hmm} that it is usually beneficial for
the first few epochs not to perform the M step in (\ref{eq:onlineEMdistrM-step})
and allow a burn-in period for the running averages of the summary
statistics to converge.

\begin{algorithm}
\selectlanguage{french}%
\begin{algorithmic}[1]

\selectlanguage{english}%
\STATE\foreignlanguage{french}{\textbf{begin}} 

\STATE\foreignlanguage{french}{\textbf{for}} $n\geq1$: let the current
parameter estimate be $\theta_{n}$. Upon obtaining measurements $Y_{n}=\{Y_{n}^{v}\}_{v\in\mathcal{V}}$
the following filtering and parameter update steps are to be performed.

\STATE\textbf{Filtering step}: Perform steps (3-6) in Algorithm \ref{alg:distFilter}.
Also compute 
\[
\tilde{\Sigma}_{n}^{r}=\left(\Sigma_{n-1}^{r}+A_{n}^{\text{T}}Q_{n}^{-1}A_{n}\right)^{-1}.
\]

\STATE\textbf{Parameter update}: Each node $r\in\mathcal{V}$ of
the network will update the following quantities for every edge $(r,j)\in\mathcal{E}$:
\begin{eqnarray*}
H_{n}^{r,j} & = & \gamma_{n}^{r}(m_{n,K}^{j,r})^{\text{T}}+(1-\gamma_{n}^{r})H_{n-1}^{r,j}\left(\tilde{\Sigma}_{n}^{r}\right)^{-1}A_{n}^{\text{T}}Q_{n}^{-1},\\
h_{n}^{r,j} & = & (1-\gamma_{n}^{r})\left(H_{n-1}^{r,j}\left(\tilde{\Sigma}_{n}^{r}\right)^{-1}\left(\Sigma_{n-1}^{r}\right)^{-1}\mu_{n-1}^{r}+h_{n-1}^{r,j}\right),\\
\mathcal{S}_{n,1}^{r,j} & = & H_{n}^{r,j}\mu_{n}^{r}+h_{n}^{r,j}.\\
\mathcal{S}_{n,2}^{r,j} & = & \gamma_{n}^{r}m_{n,K}^{j,r}+\left(1-\gamma_{n}^{r}\right)\mathcal{S}_{n-1,2}^{r,j},\\
\mathcal{S}_{n,3}^{r,j} & = & \gamma_{n}^{r}(\dot{m}_{n,K}^{j,r}-\ddot{m}_{n,K}^{j,r}+\ddot{m}_{n,1}^{j,r})+\left(1-\gamma_{n}^{r}\right)\mathcal{S}_{n-1,3}^{r,j},
\end{eqnarray*}

Upon doing so the localization parameter is updated:
\begin{equation}
\theta_{n+1}^{r,j}=\Lambda(\mathcal{S}_{n,1}^{r,j},\mathcal{S}_{n,2}^{r,j},\mathcal{S}_{n,3}^{r,j}).\label{eq:onlineEMdistrM-step}
\end{equation}

\STATE\foreignlanguage{french}{\textbf{end for}}

\STATE\foreignlanguage{french}{\textbf{end}}

\selectlanguage{french}%
\end{algorithmic}

\selectlanguage{english}%
\caption{\textbf{Distributed on-line EM}}
\label{alg:onlineEM}
\end{algorithm}

\section{Numerical Examples\label{sec:examples}}

\begin{figure}
\centering

\subfigure[11 node sensor network]{ \includegraphics[width=0.3\textwidth]{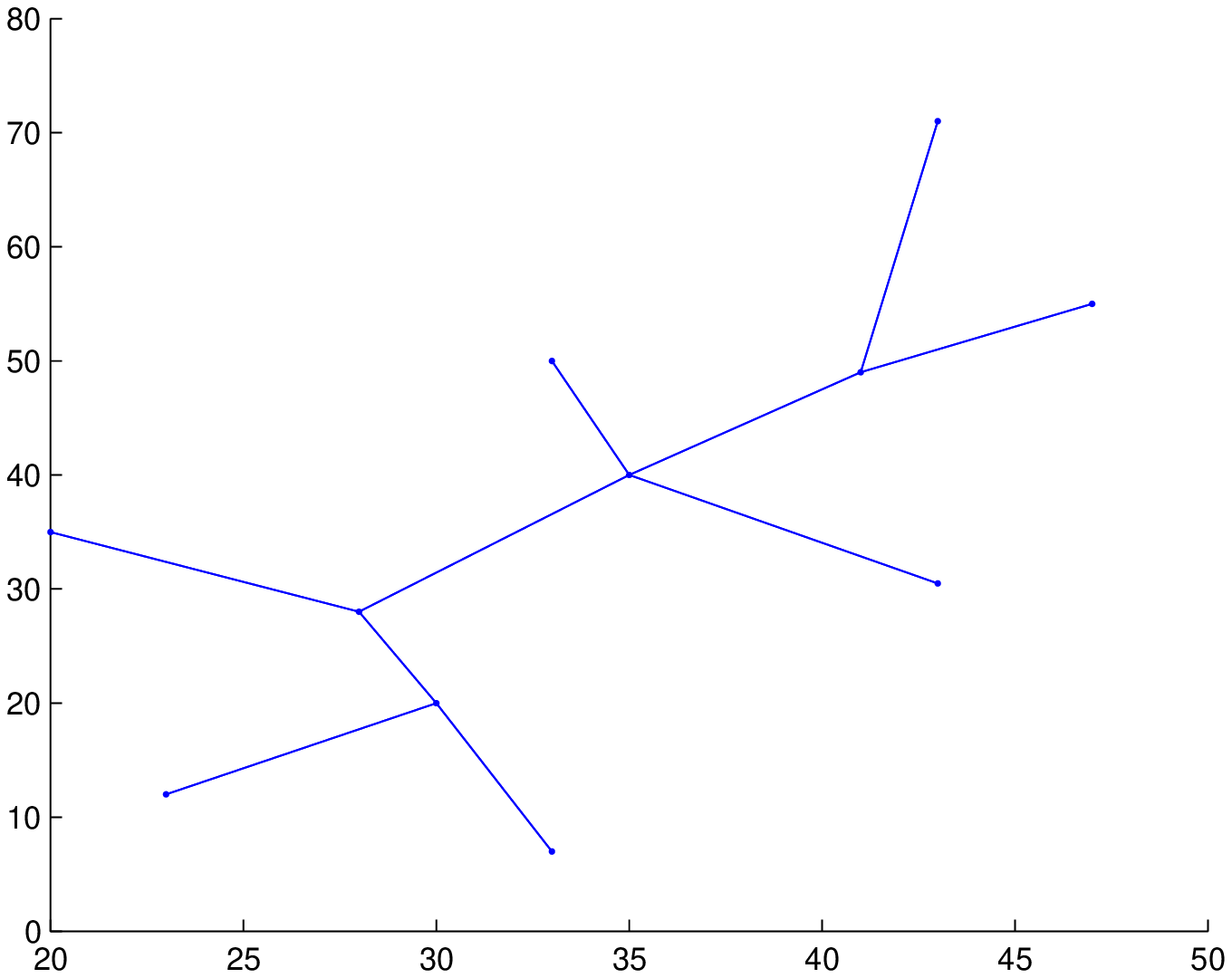}
\label{fig:sn1}}\subfigure[44 node sensor network]{ \includegraphics[width=0.3\textwidth]{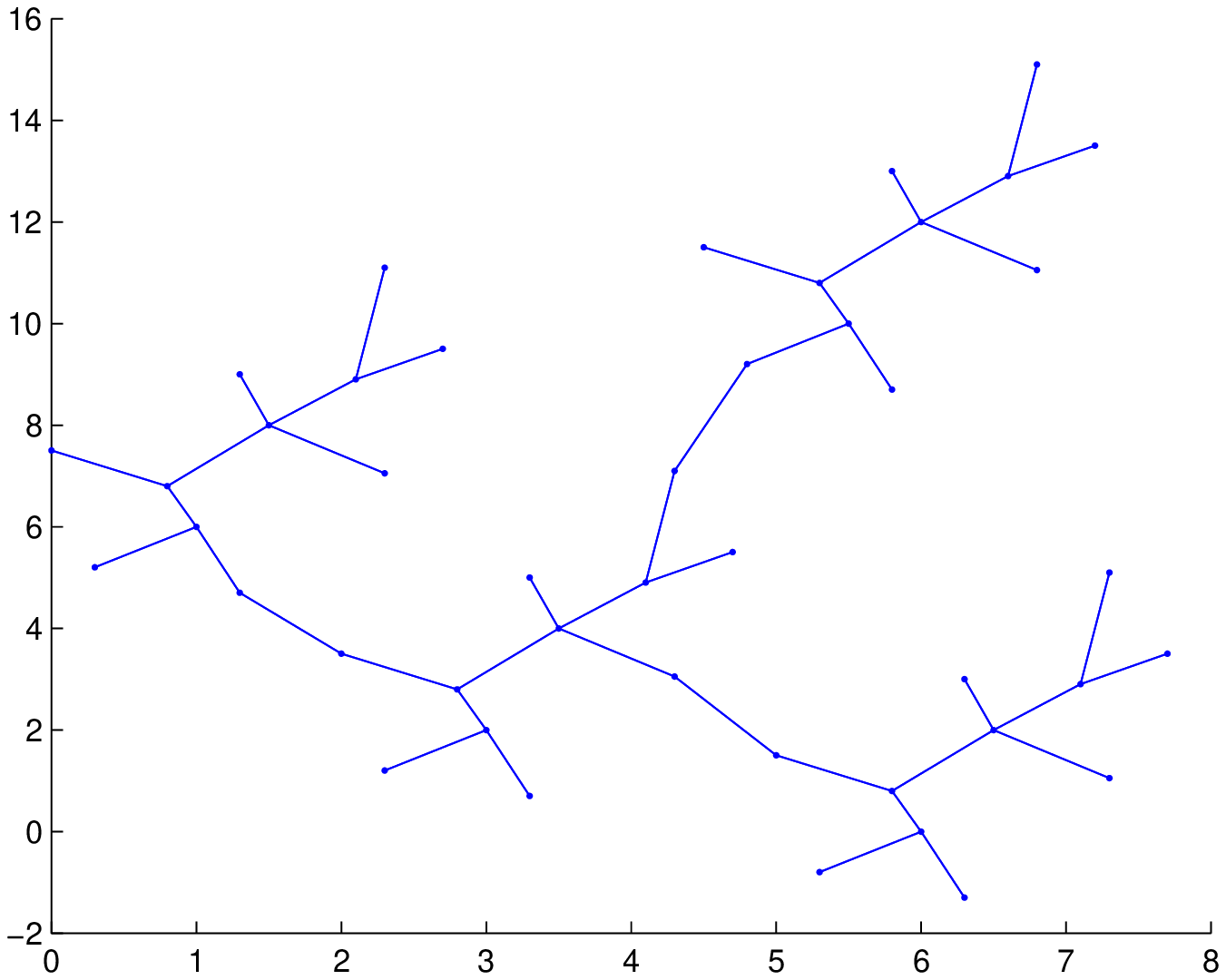}
\label{fig:sn2} }\subfigure[11 node sensor network with cycles]{\includegraphics[width=0.3\textwidth]{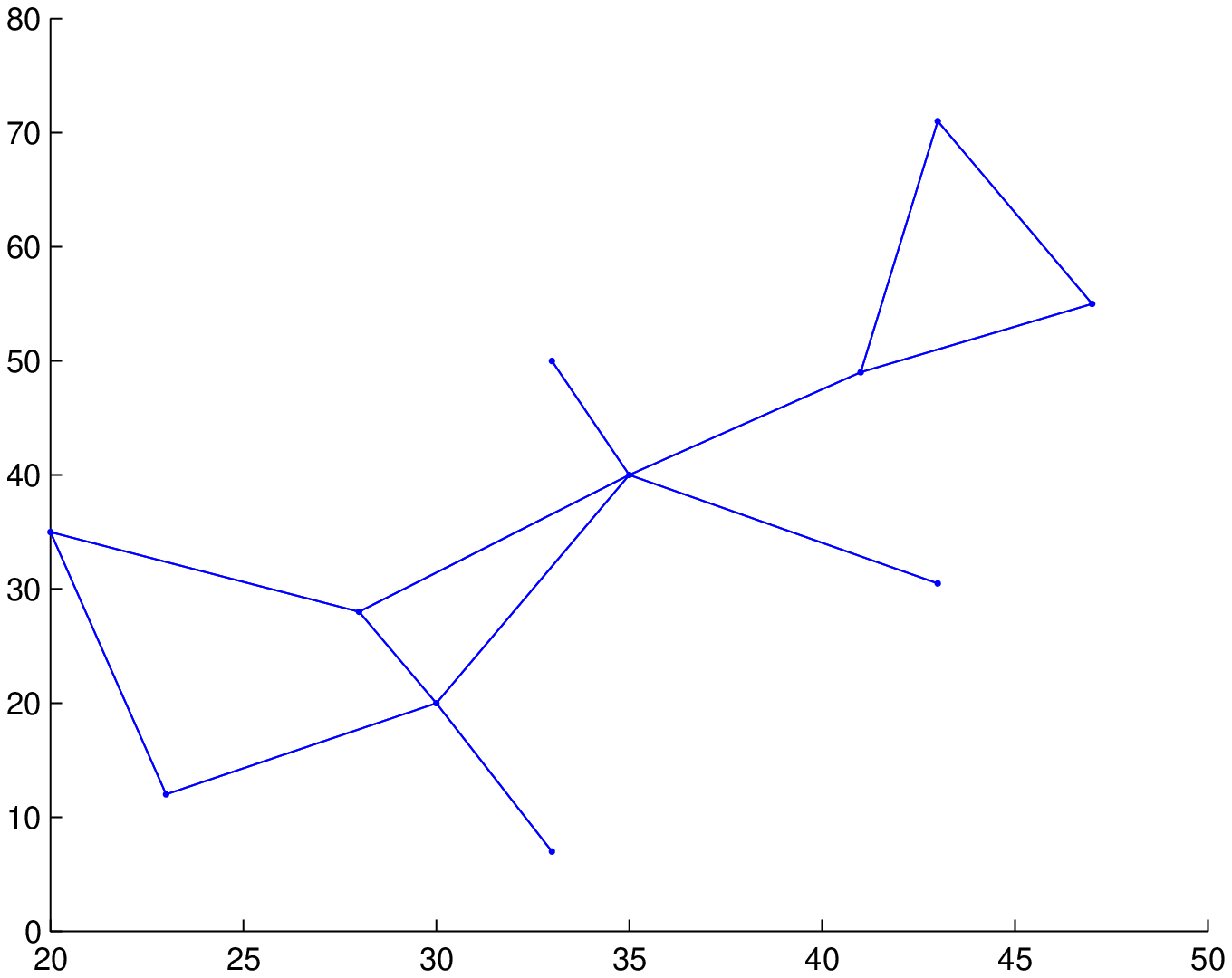}\label{fig:loopysn}}

\caption{Various sensor networks of different size and topology.}

\label{fig:sub1} 
\end{figure}

The performance of the distributed RML and EM algorithms are studied
using a Linear Gaussian and a non-linear model. For both cases the
hidden target is given in (\ref{eq:linearGaussianState}) with $V_{n}=B\widetilde{V}_{n},$\ where
$\widetilde{V}_{n}$ is zero mean Gaussian additive noise with variance
$\widetilde{Q}_{n}$, and 
\[
A_{n}=\left[\begin{array}{cccc}
1 & \tau & 0 & 0\\
0 & 1 & 0 & 0\\
0 & 0 & 1 & \tau\\
0 & 0 & 0 & 1
\end{array}\right],\quad B=\left[\begin{array}{cc}
\frac{\tau^{2}}{2} & 0\\
\tau & 0\\
0 & \frac{\tau^{2}}{2}\\
0 & \tau
\end{array}\right],\quad\widetilde{Q}_{n}=\sigma_{x}^{2}I,
\]
 and $I$ is the identity matrix. %
{} For the linear model the observations are given by (\ref{eq:linearGaussianObs})
with 
\[
C_{n}^{r}=\alpha^{r}\left[\begin{array}{cccc}
1 & 0 & 0 & 0\\
0 & 0 & 1 & 0
\end{array}\right],\quad R_{n}^{r}=\sigma_{y}^{2}I,
\]
 where $\alpha^{r}$ are constants different for each node and are
assigned randomly from the interval $[0.75,1.25]$. For the non-linear
model we will use the bearings-only measurement model. In this model
at each node $r$, the observation $Y_{n}^{r}$\ is: 
\[
Y_{n}^{r}=tan^{-1}(X_{n}^{r}(1)/X_{n}^{r}(3))+W_{n}^{r}.
\]
 with $W_{n}^{r}\overset{i.i.d.}{\sim}\mathcal{N}(0,{0.35}^{2})$.
For the remaining parameters we set $\tau=0.01$, $\sigma_{x}$ $=$1
and $\theta_{0}^{r,j}=0$ for all $(r,j)\in\mathcal{E}$ . In Figure
\ref{fig:sub1} we show three different sensor networks for which
we will perform numerical experiments.

\begin{figure}
\centering

\subfigure[{RML for tree network}]{ \includegraphics[width=0.3\textwidth]{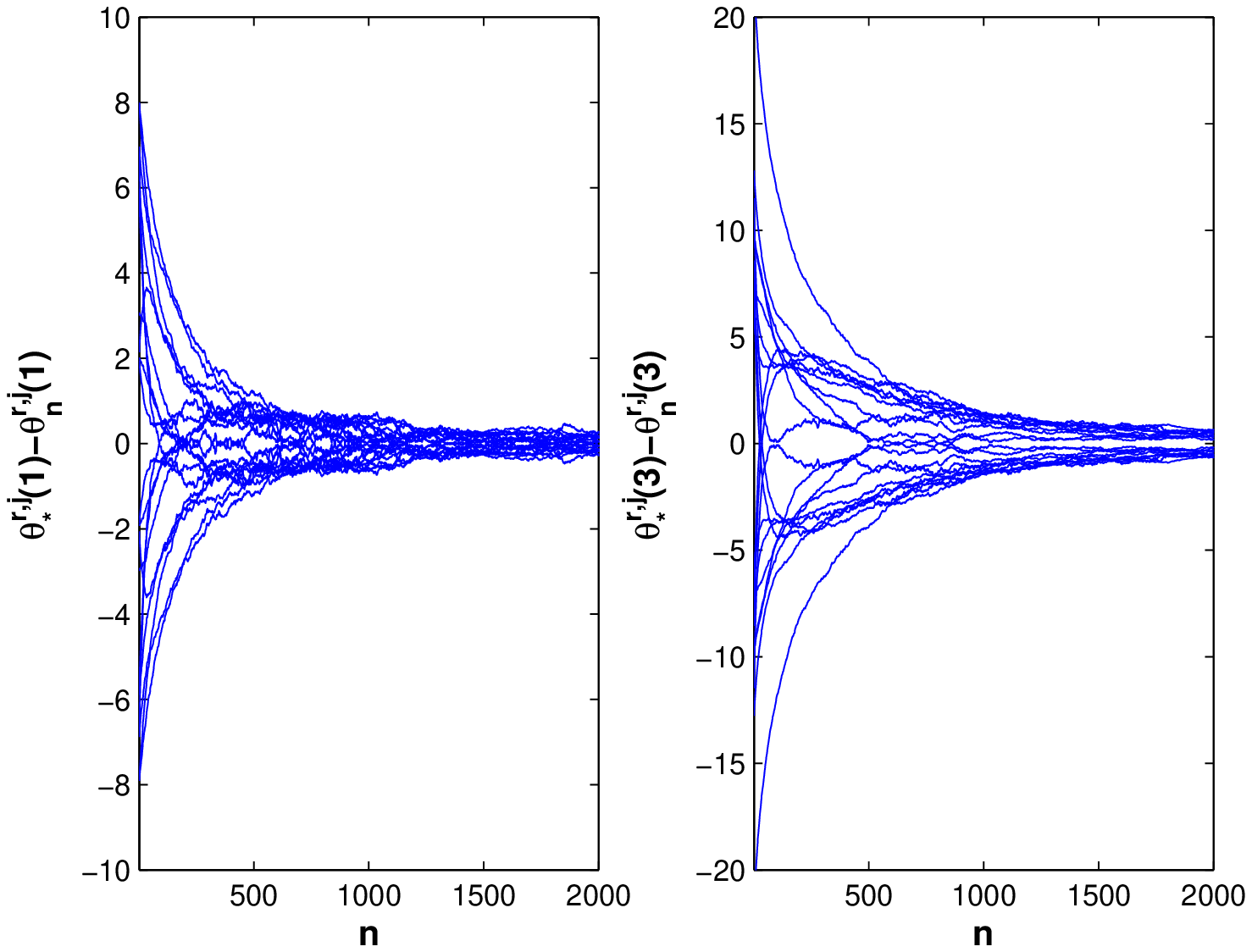}
\label{fig:rml} }\subfigure[{Nonlinear RML for tree network}]{
\includegraphics[width=0.3\textwidth]{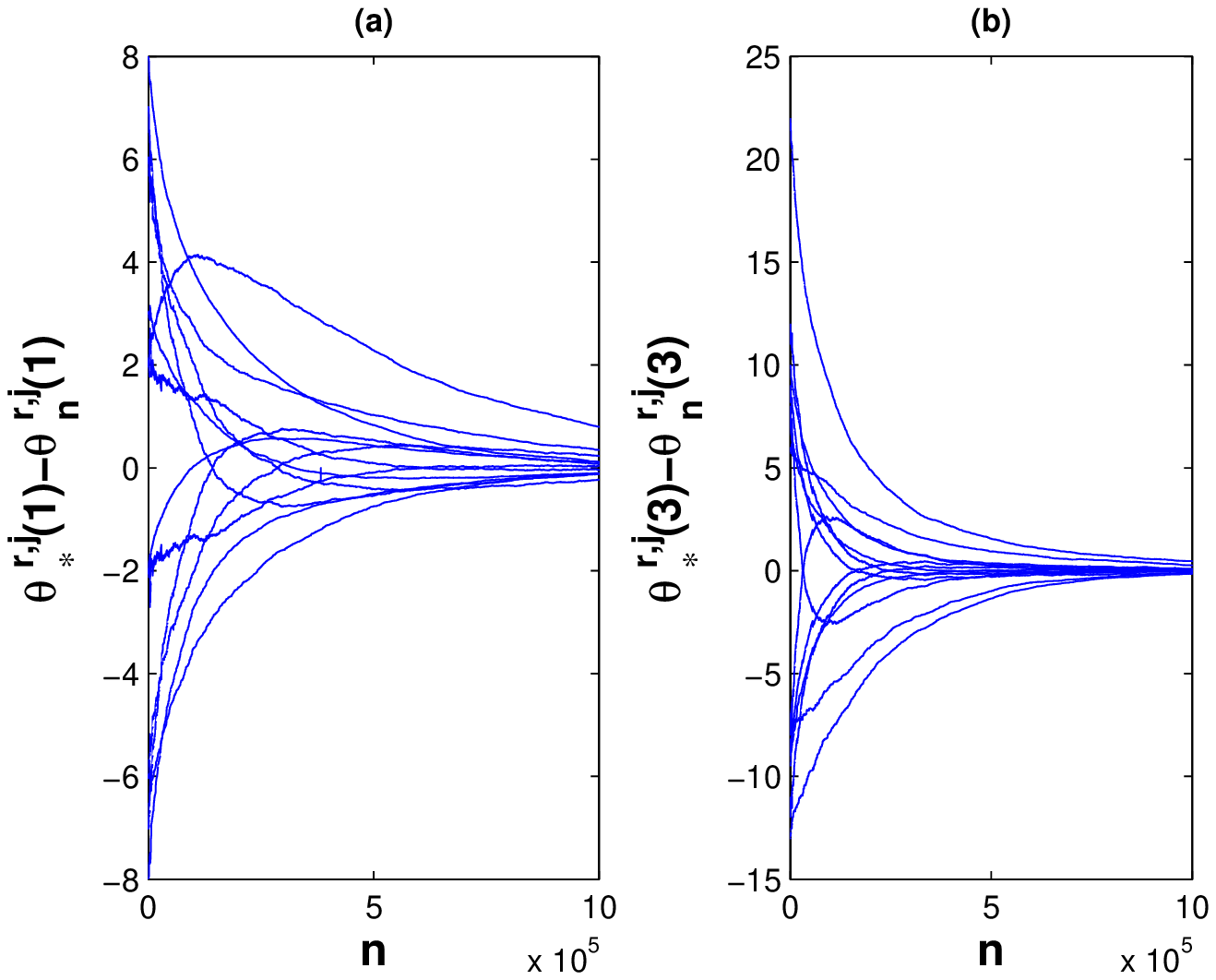} \label{fig:ekf}
}\subfigure[{RML for network with cycles}]{\includegraphics[width=0.3\textwidth]{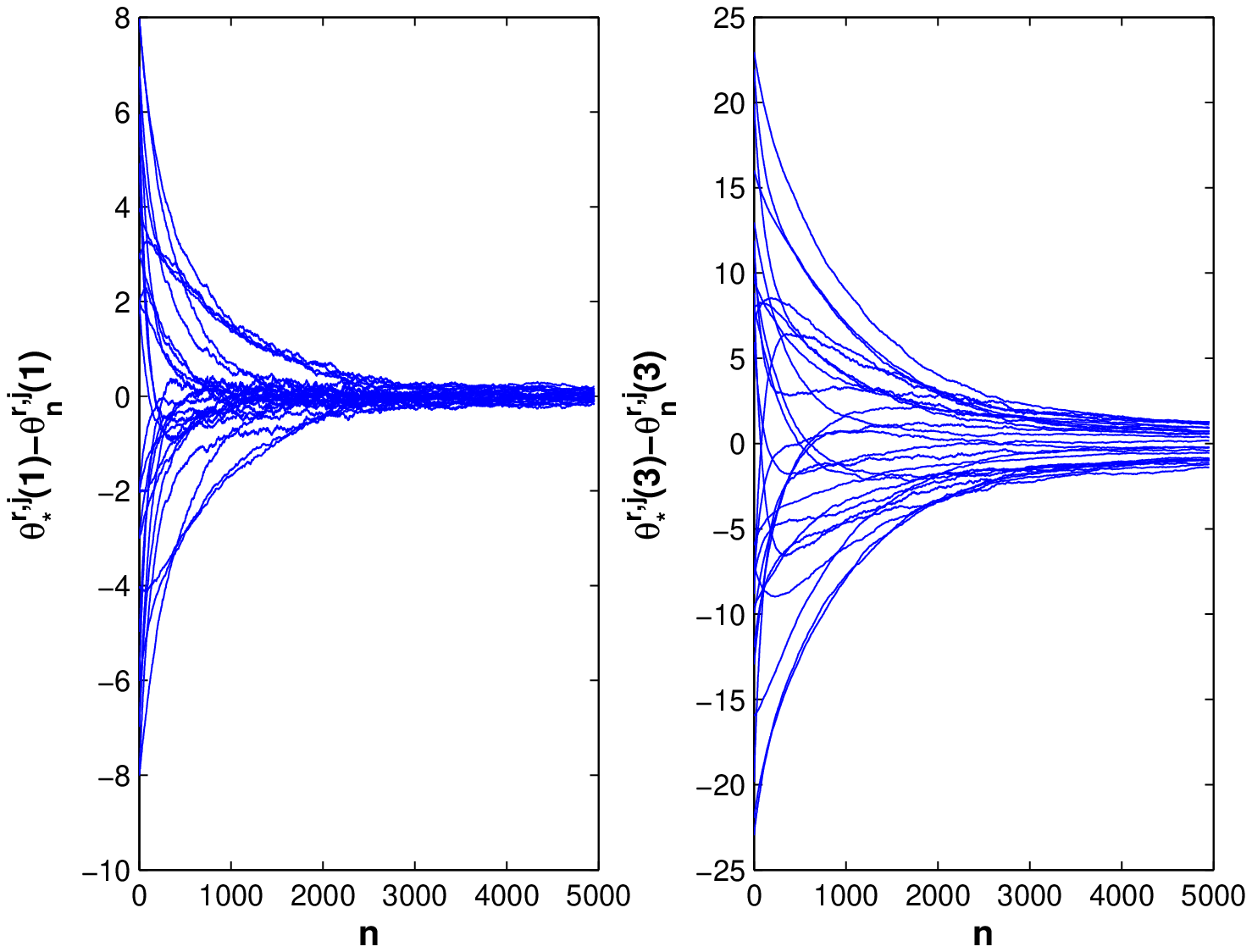}\label{fig:rmlloop}}\\

\subfigure[{EM for tree network}]{ \includegraphics[width=0.3\textwidth]{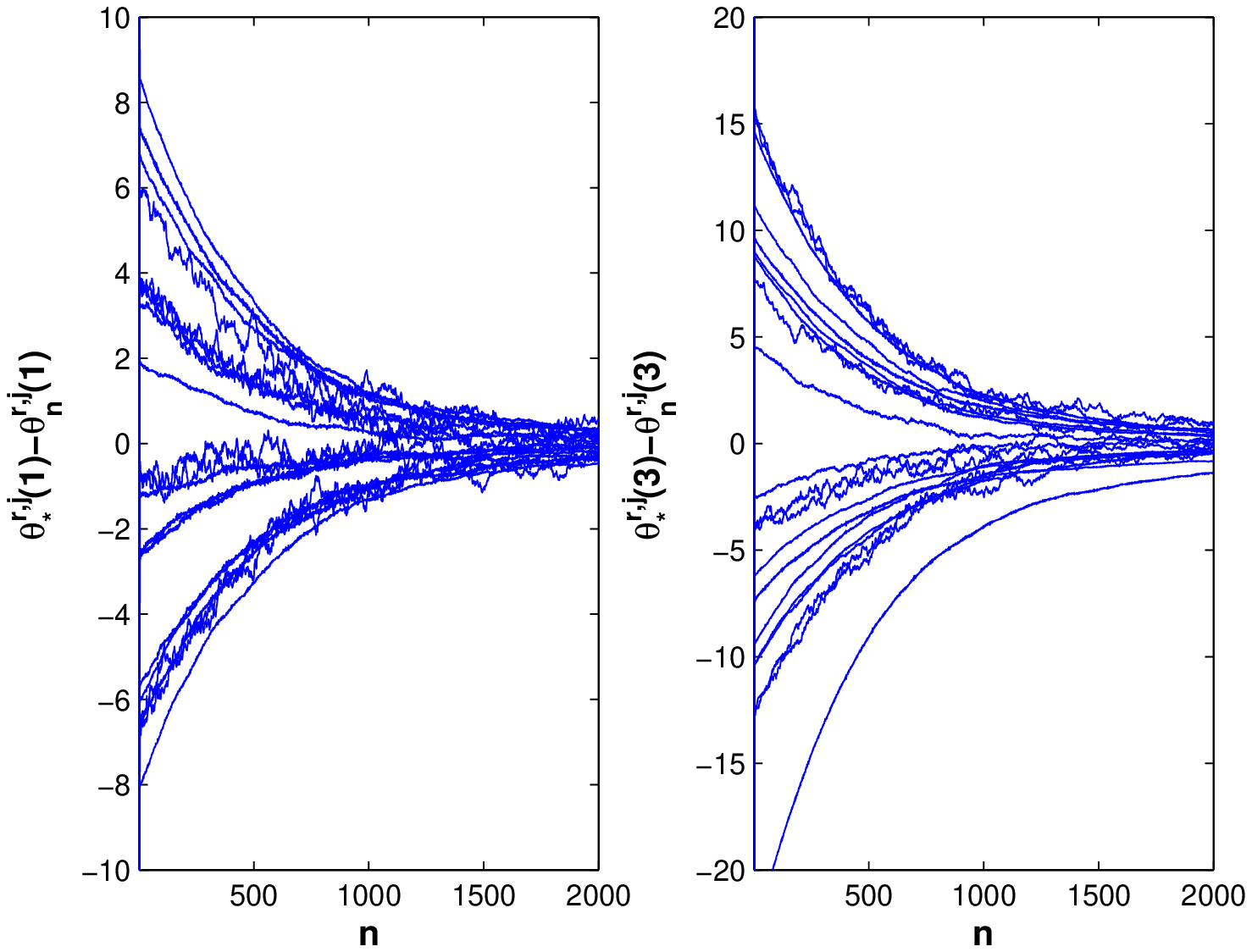}
\label{fig:em} } \subfigure[{Nonlinear RML for large network}]{
\includegraphics[width=0.3\textwidth]{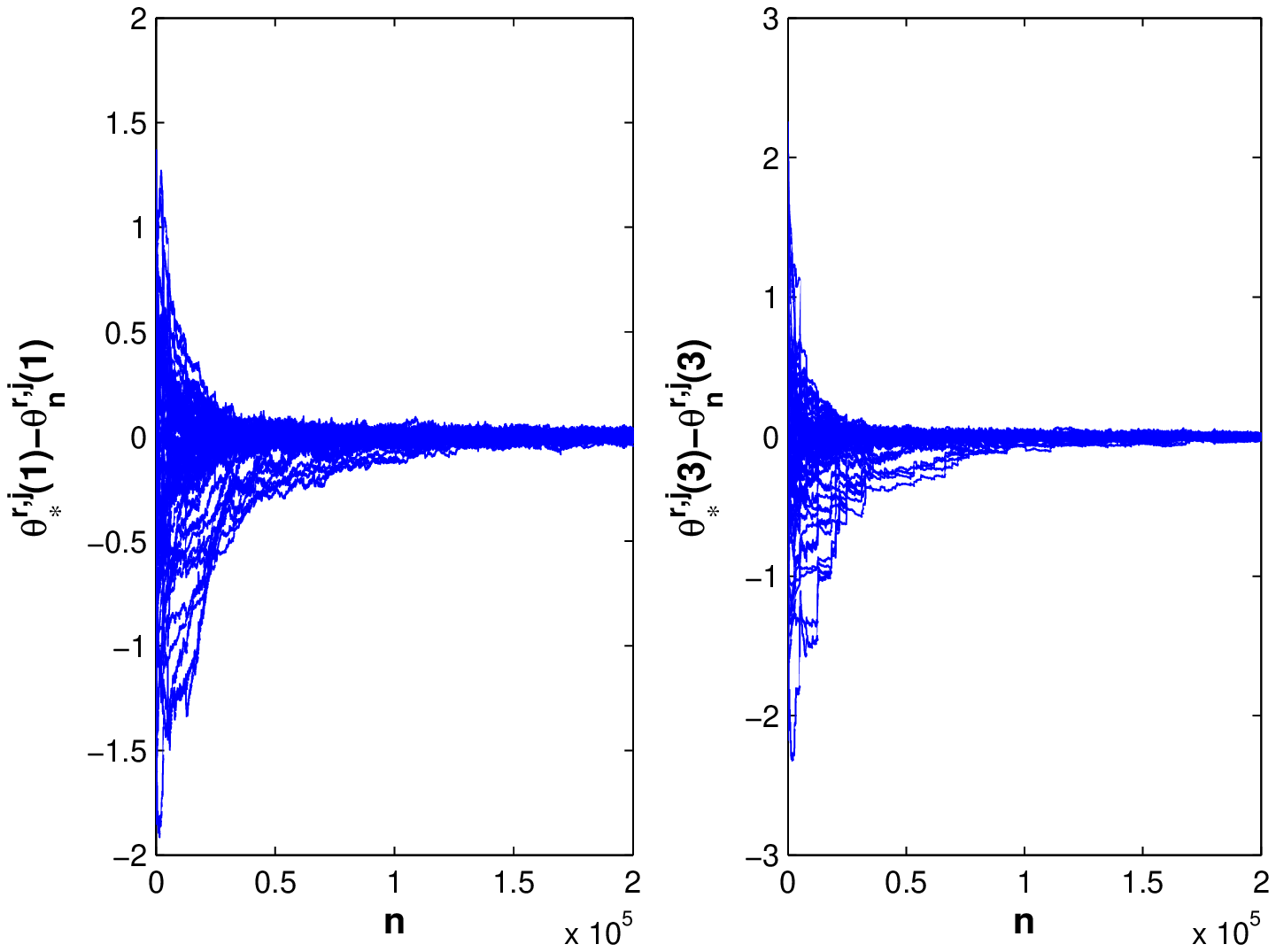} \label{fig:ekf2}
}\subfigure[{EM for network with cycles}]{\includegraphics[width=0.3\textwidth]{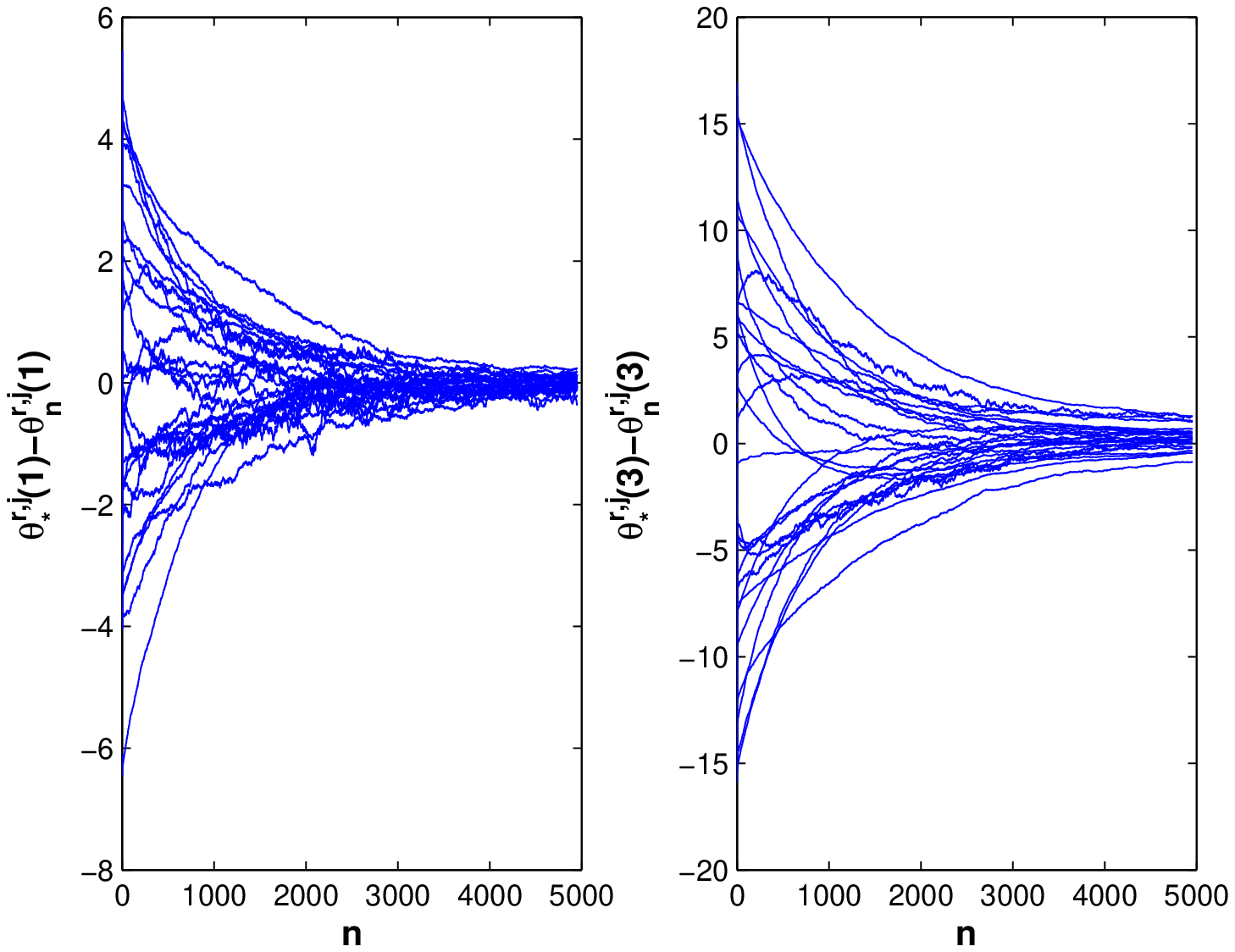}\label{fig:emloop}}

\caption{The convergence of the localization parameters' estimates to $\theta_{\ast}^{r,j}$
is demonstrated using appropriate error plots for various sensor networks.
Left: Parameter error after each iteration for each edge of the medium
sensor network of Fig. \ref{fig:sn1}. In each subfigure left and
right columns show the errors in the x- and y- coordinates respectively;
(a) is for RML and (d) is for EM. Middle: Same errors when using RML
for the nonlinear bearings-only observation model; (b) is for medium
sized network of Fig. \ref{fig:sn1} and (e) for the large network
of Fig. \ref{fig:sn2}. Right: Same errors for network with cycles
seen in Fig \ref{fig:loopysn}; (c) for RML and (f) for EM. }

\label{fig:sub2} 
\end{figure}

In Figure \ref{fig:sub2} we present various convergence plots for
each of these networks for $\sigma_{y}=0.5$. We plot both dimensions
of the errors $\theta_{\ast}^{r,j}-\theta_{n}^{r,j}$ for three cases:
\begin{itemize}
\item in (a) and (d) we use distributed RML and on-line EM respectively
for the network of Figure \ref{fig:sn1} and the linear Gaussian model. 
\item in (b) and (e) we use distributed RML for the bearings only tracking
model and the networks of Figures \ref{fig:sn1} and \ref{fig:sn2}
respectively. Local linearization as discussed in Sections \ref{sub:Non-linear-Model},
\ref{sec:distFiltering} and \ref{sec:distRML} was used to implement
the distributed RML algorithm. We remark that we do not apply the
online EM to problems where the solution to the M-step cannot be expressed
analytically as some function $\Lambda$ of summary statistics. %

\item in (c) and (f) we use distributed RML and on-line EM for respectively
for the network of Figure \ref{fig:loopysn} and the linear Gaussian
model. In this case we used $K=2$. 
\end{itemize}
All errors converge to zero. Although both methods are theoretically
locally optimal when performing the simulations we did not observe
significant discrepancies in the errors for different initializations.
For both RML and on-line EM we used for $n\leq10^{3}$ a constant
but small step-size, $\gamma_{n}^{r}=\gamma=4\times10^{-3}$ and $0.025$
respectively. For the subsequent iterations we set $\gamma_{n}^{r}=\gamma(n-10^{3}){}^{-0.8}$
. Note that if the step-size decreases too quickly in the first time
steps, these algorithms might converge too slowly. In the plots of
Figure \ref{fig:sub2} one can notice that the distributed RML and
EM algorithms require comparable amount of time to converge with the
RML being usually faster. For example in Figures \ref{fig:sub2} (a)
and (d) we observe that RML requires around $1000$ iterations to
converge whereas on-line EM requires approximately 2000 iterations.
We note that the converge rate also depends on the specific network
used, the value of $K$ and the simulation parameters. 

\begin{figure}
\subfigure[{RML for different K }]{\includegraphics[width=0.25\textwidth]{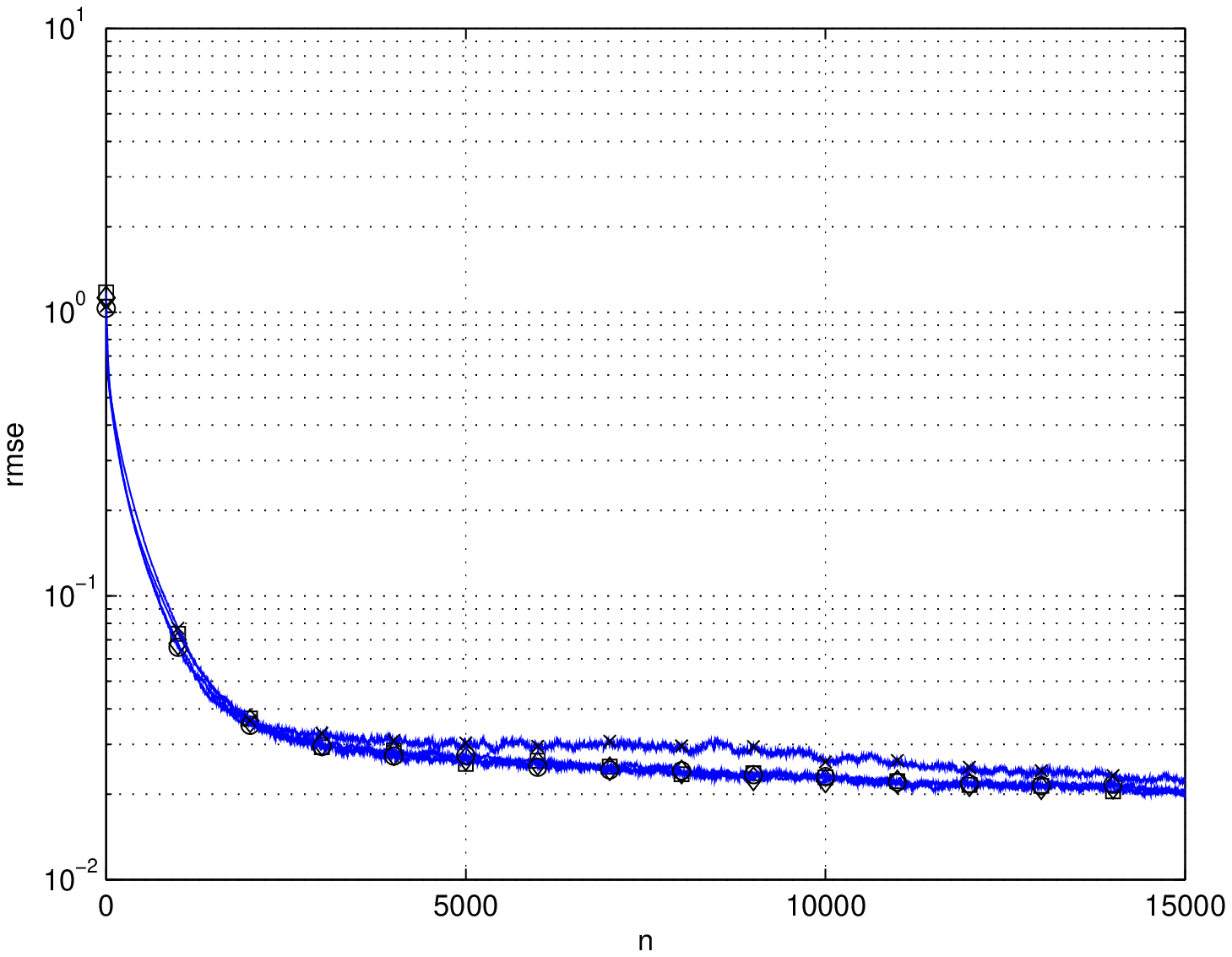}}\subfigure[{EM for different K }]{\includegraphics[width=0.25\textwidth]{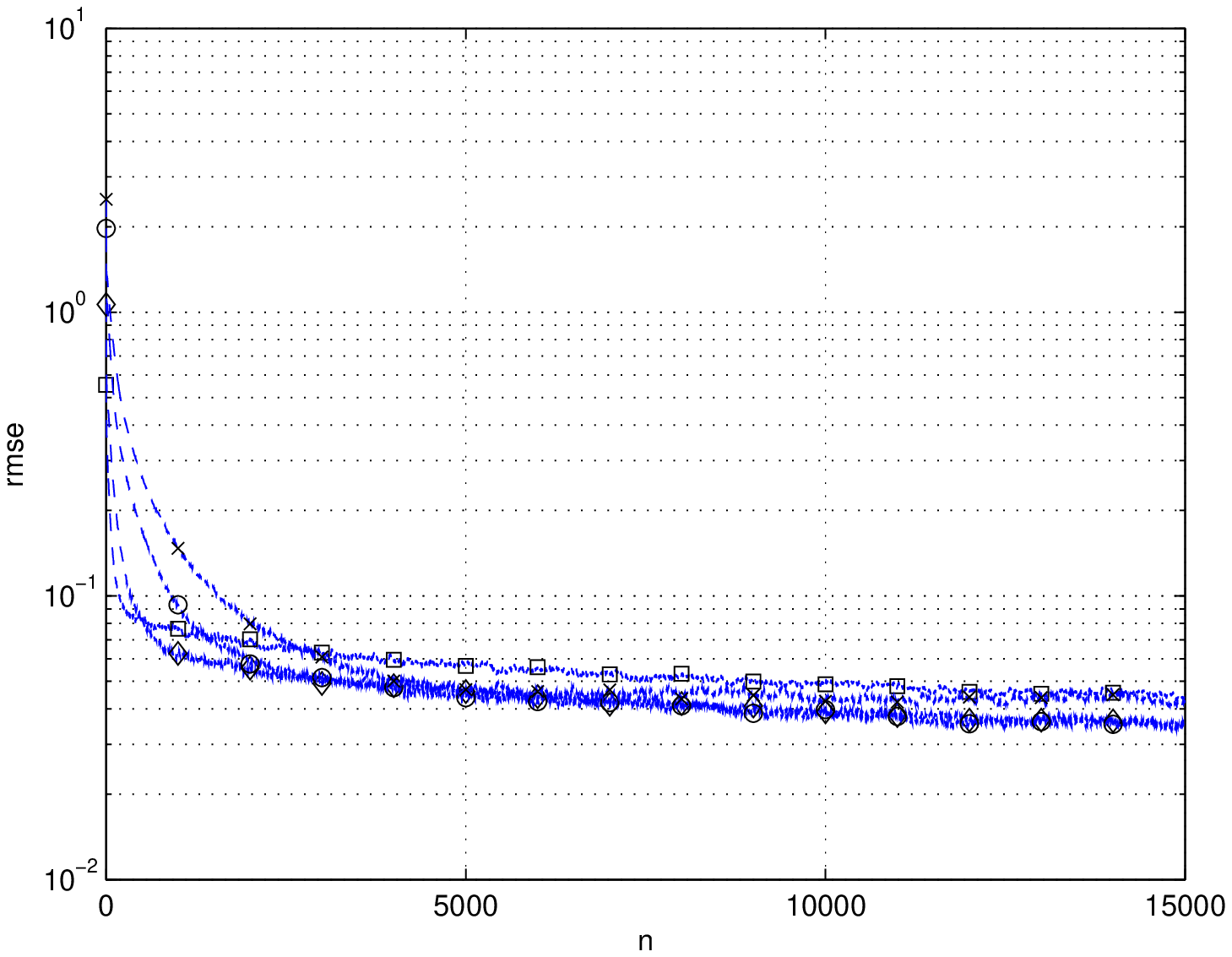}}\subfigure[{RML for different $\frac{\sigma_{x}}{\sigma_{y}}$ }]{\includegraphics[width=0.25\textwidth]{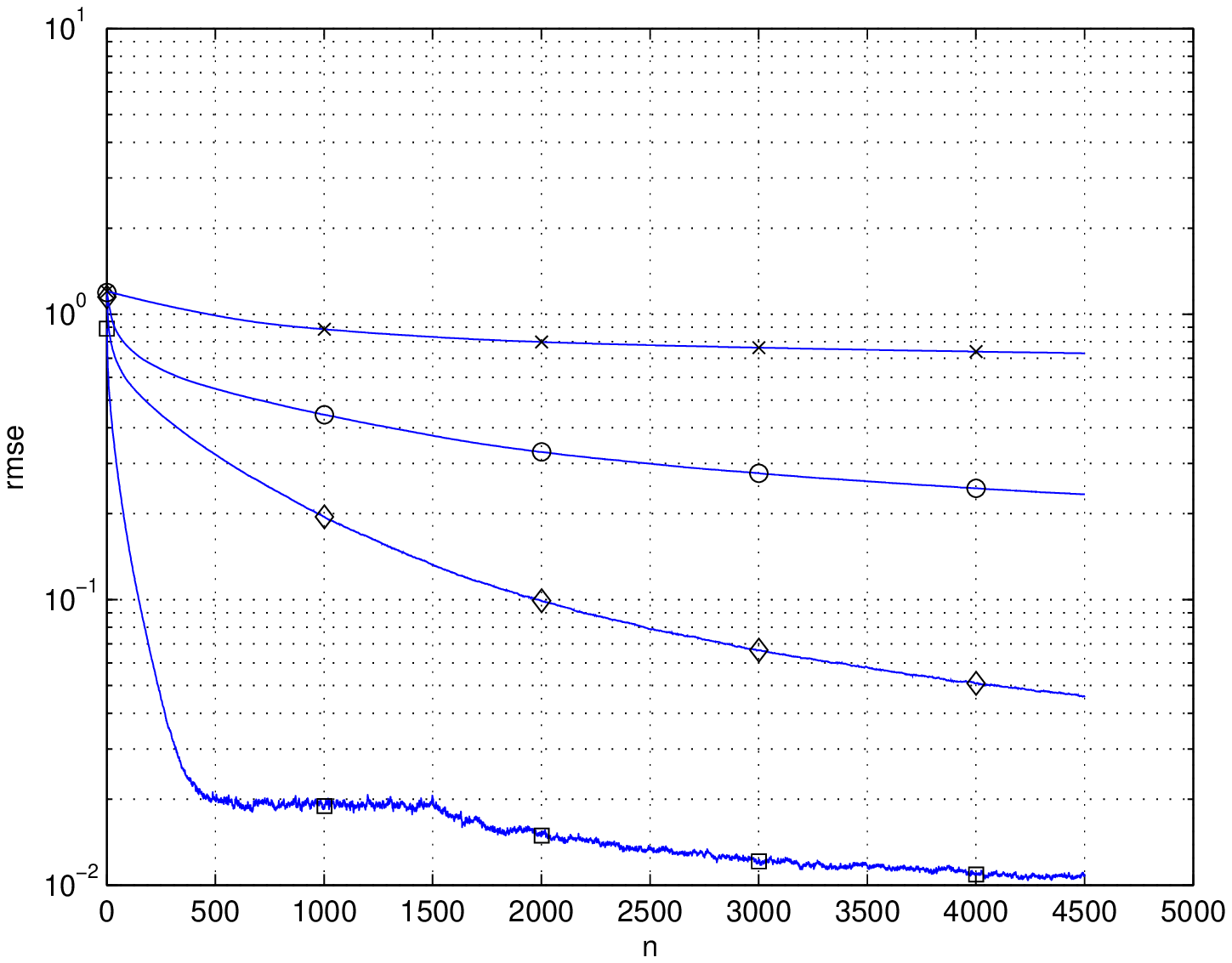}}\subfigure[{EM for different $\frac{\sigma_{x}}{\sigma_{y}}$ }]{\includegraphics[width=0.25\textwidth]{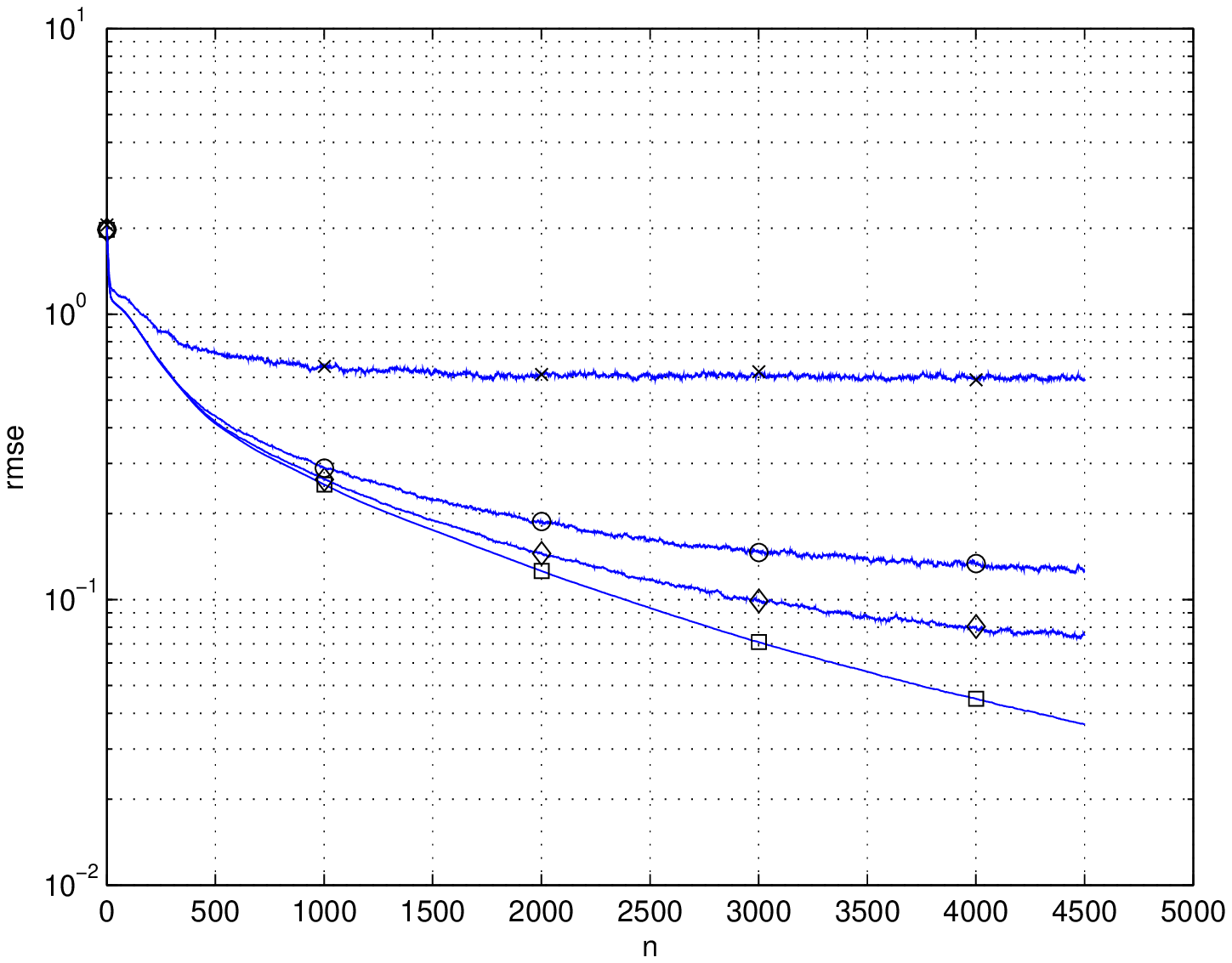}}\caption{Comparison of distributed RML and on-line EM. (a) (and (b) resp.):
RMSE for RML (and on-line EM resp.)$ $ against $n$ for $K=$ $2$
($\square$), $4$ ($\diamond$), $8$ ($\circ$), 12($\times$).
(c) (and (d) resp.): RMSE for RML (and on-line EM resp.)$ $ for $\frac{\sigma_{x}}{\sigma_{y}}=$
$10$ ($\square$), $1$ ($\diamond$), $0.5$ ($\circ$), $0.1$($\times$). }
\label{sub:snr_K}
\end{figure}

To investigate this further we varied $K$ and $\frac{\sigma_{x}}{\sigma_{y}}$
and recorded the root mean squared error (RMSE) for $\theta_{n}$
obtained for the network of Figure \ref{fig:sn2} using $50$ independent
runs. For the RMSE at time $n$ we will use $\sqrt{\frac{1}{50\left|\mathcal{E}\right|}\sum_{e\in\mathcal{E}}\sum_{m=1}^{50}\left\Vert \theta_{\ast}^{r,j}-\theta_{n,m}^{r,j}\right\Vert _{2}^{2}}$,
where $\theta_{n,m}^{r,j}$ denotes the estimated parameter at epoch
$n$ obtained from the $m$-th run. The results are plotted in Figure
\ref{sub:snr_K} for different cases:
\begin{itemize}
\item in (a) and (b) for $\frac{\sigma_{x}}{\sigma_{y}}=2$ we show the
RMSE for $K=2,4,8,12$. We observe that in every case the RMSE keeps
reducing as $n$ increases. Both algorithms behave similarly with
the RML performing better and showing quicker convergence. One expects
that observations beyond your near immediate neighbors are not necessary
to localize adjacent nodes and hence the good performance for small
values of $K$.
\item in (b) and (c) we show the RMSE for RML and on-line EM respectively
when $\frac{\sigma_{x}}{\sigma_{y}}=10,1,0.5,0.1$. We observe that
EM seems to be slightly more accurate for lower values of $\frac{\sigma_{x}}{\sigma_{y}}$
with the reverse holding for higher values of the ratio.
\end{itemize}
In each run the same step-size was used as before except for RML and
$\frac{\sigma_{x}}{\sigma_{y}}=10$, where we had to reduce the step
size by a factor of $10$. 

\section{Conclusion\label{sec:Discussion}}

We have presented a method to perform collaborative tracking and self-localization.
We exploited the fact that different nodes collect measurements of
a common target. This idea has appeared previously in \cite{Funiak06,TaylorLaSlat},
both of which use a Bayesian inference scheme for the localization
parameters. We remark that our distributed ML methods appear simpler
to implement than these Bayesian schemes as the messages here are
nothing more than the appropriate summary statistics for computing
the filtering density and performing parameter updates. There is good
empirical evidence that the distributed implementations of ML proposed
in this paper are stable and do seem to settle at reasonably accurate
estimates. %
{} A theoretical investigation of the properties of the schemes would
be an interesting but challenging extension. %
{} Finally, as pointed out by one referee, another interesting extension
would be to develop consensus versions of Algorithm 1 in the spirit
of gossip algorithms in \cite{kempe2003gossip} or the aggregation
algorithm of \cite{bertsekas1989parallel} %
{} which might be particularly relevant for networks with cycles, which
are dealt with here by using an appropriate value for $K$.


\section{Maximum likelihood parameter estimation\label{sub:MLE}}

This section does not pertain to sensor localization specifically
but to the general problem of static parameter estimation in HMMs
using ML. Thus to avoid confusion with the localization problem a
different font is used the notation. Consider a HMM where $\{\mathsf{X}_{n}\}_{n\geq1}$\ is
the hidden state-process and $\{\mathsf{Y}_{n}\}_{n\geq1}$\ is the
observed process each taking values in taking values in $\mathbb{R}^{d_{x}}$
and $\mathbb{R}^{d_{y}}$ respectively. For the transition density
for $\{\mathsf{X}_{n}\}_{n\geq1}$, we have $\mathsf{X}_{n+1}|\mathsf{X}_{n}=\mathsf{x}_{n}\sim\mathsf{f}(\cdot|\mathsf{x}_{n})$.
The observation model, $\mathsf{Y}_{n}|\mathsf{X}_{n}=\mathsf{x}_{n}\sim\mathsf{g}_{\vartheta}(\cdot|\mathsf{x}_{n})$
is parametrized by $\vartheta\in\Theta\:(\subset\mathbb{R}^{d_{\vartheta}})$.
The true static parameter generating the sequence of observations
is $\vartheta_{\ast}$\ and is to be learned from the observed data
$\{\mathsf{Y}_{n}\}_{n\geq1}$. The ML parameter estimate is the maximizing
argument of the log-likelihood of the observed data up to time $n$:
$\tilde{\vartheta}_{n}=\arg\max_{\vartheta\in\Theta}\log p_{\vartheta}(\mathsf{Y}_{1:n})$.
Here $p_{\vartheta}(\mathsf{Y}_{1:n})$ denotes the joint density
of $\mathsf{Y}_{1:n}$ and the subscript makes explicit the value
of the parameter used to compute this density.

For a long observation sequence we are interested in a recursive parameter
estimation procedure in which the data is run through once sequentially.
If $\vartheta_{n}$ is the estimate of the model parameter after $n$
observations, a recursive method would update the estimate to $\vartheta_{n+1}$\ after
receiving the new data $\mathsf{Y}_{n}$. For example, consider the
following update scheme: 
\begin{equation}
\vartheta_{n+1}=G_{n+1}(\vartheta_{n},\mathsf{Y}_{n}),\qquad n\geq1.\label{eq:titsRML}
\end{equation}
 where $G_{n+1}$ is an appropriate function to be defined. This scheme
was originally suggested by \cite{Tit84,TiJ83} when $\{\mathsf{X}_{n}\}_{n\geq1}$
is not a Markov chain but rather an independent and identically distributed
(i.i.d.) sequence. 


\subsection{Recursive Maximum Likelihood (RML) \label{subsec:RML}}

To motivate a suitable choice for $G_{n+1}(\vartheta_{n},\mathsf{Y}_{n})$
for estimating the parameters of a HMM, consider the following recursion:
\begin{equation}
\vartheta_{n+1}=\vartheta_{n}+\gamma_{n+1}\left.\nabla\log p_{\vartheta}(\mathsf{Y}_{n}|\mathsf{Y}_{1:n-1})\right\vert _{\vartheta=\vartheta_{n}}.\label{eq:almostOnlineRML}
\end{equation}
 where $\{\gamma_{n}\}$ is the step-size sequence that should satisfy
the following constraints: $\sum_{n}\gamma_{n}=\infty$ and $\sum_{n}\gamma_{n}^{2}<\infty$.
One possible choice would be $\gamma_{n}=n^{-\alpha}$, $0.5<\alpha<1$.
Here $p_{\vartheta}(\mathsf{Y}_{n}|\mathsf{Y}_{1:n-1})$ is the conditional
density of $\mathsf{Y}_{n}$\ given $\mathsf{Y}_{1:n-1}\ $ and the
subscript makes explicit the value of the parameter used to compute
this density. Upon receiving $\mathsf{Y}_{n}$, $\vartheta_{n}$\ is
updated in the direction of ascent of the conditional density of this
new observation. The algorithm in the present form is not suitable
for online implementation due to the need to evaluate the gradient
of $\log p_{\vartheta}(\mathsf{Y}_{n}|\mathsf{Y}_{1:n-1})$ (w.r.t.
$\vartheta$)\ at $\vartheta=\vartheta_{n}$. Doing so would require
browsing through the entire history of observations. This limitation
is removed by defining certain intermediate quantities that facilitate
the online evaluation of this gradient \cite{LeGlandMevel97,CoR98}.

In particular, assume that from the previous iteration of the RML,
one has computed $p_{n}(\mathsf{x}_{n})\approx\left.p_{\vartheta}(\mathsf{x}_{n}|\mathsf{Y}_{1:n-1})\right\vert _{\vartheta=\vartheta_{n}}$
and $\dot{p}_{n}(\mathsf{x}_{n})\approx\left.\nabla p_{\vartheta}(\mathsf{x}_{n}|\mathsf{Y}_{1:n-1})\right\vert _{\vartheta=\vartheta_{n}.}$,
where $(p_{n},\dot{p}_{n})$ are approximations of the predicted density
and its gradient evaluated at $\vartheta=\vartheta_{n}$. The RML
is initialized with an arbitrary value for $\vartheta_{1}$, $p_{1}(\mathsf{x}_{1})=p_{\vartheta_{1}}(\mathsf{x}_{1})$,
which is the prior distribution for $\mathsf{X}_{1}$ and $\dot{p}_{1}(\mathsf{x}_{1})=\left.\nabla p_{\vartheta}(\mathsf{x}_{1})\right\vert _{\vartheta=\vartheta_{1}}$,
i.e. the gradient of this prior which could be zero if it does not
depend on $\vartheta$. Then the online version of (\ref{eq:almostOnlineRML}),
which is the RML procedure of \cite{LeGlandMevel97,CoR98}, proceeds
as follows. Given the new observation $\mathsf{Y}_{n}$, update the
parameter: 
\begin{equation}
\vartheta_{n+1}=\vartheta_{n}+\gamma_{n+1}\left(\int\mathsf{g}_{\vartheta_{n}}(\mathsf{Y}_{n}|\mathsf{x}_{n})p_{n}(\mathsf{x}_{n})d\mathsf{x}_{n}\right)^{-1}\left(\int\mathsf{\dot{g}}_{\vartheta_{n}}(\mathsf{Y}_{n}|\mathsf{x}_{n})p_{n}(\mathsf{x}_{n})d\mathsf{x}_{n}+\int\mathsf{g}_{\vartheta_{n}}(\mathsf{Y}_{n}|\mathsf{x}_{n})\dot{p}_{n}(\mathsf{x}_{n})d\mathsf{x}_{n}\right)\label{eq:centralrml}
\end{equation}
 where $n\geq1$ and $\mathsf{\dot{g}}_{\vartheta^{\prime}}(\mathsf{y}|\mathsf{x})\equiv\left.\nabla_{\vartheta}\mathsf{g}_{\vartheta}(\mathsf{y}|\mathsf{x})\right\vert _{\vartheta=\vartheta^{\prime}}$.
In (\ref{eq:almostOnlineRML}), the desired gradient is the ratio
of the terms $\left.p_{\vartheta}(\mathsf{Y}_{n}|\mathsf{Y}_{1:n-1})\right\vert _{\vartheta=\vartheta_{n}}$
and $\left.\nabla p_{\vartheta}(\mathsf{Y}_{n}|\mathsf{Y}_{1:n-1})\right\vert _{\vartheta=\vartheta_{n}}$.
This ratio is approximated in the fraction on the right-hand side
of (\ref{eq:centralrml}). After computing (\ref{eq:centralrml}),
one may update $(p_{n},\dot{p}_{n})$ to $(p_{n+1},\dot{p}_{n+1})$
for the next RML iteration. Specific expressions for this update may
be found for example in \cite[Section 8.2.1]{kantas-thesis} or \cite{LeGlandMevel97}.
The recursive propagation of $(p_{n},\dot{p}_{n})$ implicitly involves
the previous values of the parameter, i.e. $\vartheta_{1:n}$, and
hence are only approximations to $\left.p_{\vartheta}(\mathsf{x}_{n}|\mathsf{Y}_{1:n-1})\right\vert _{\vartheta=\vartheta_{n+1}}$,
$\left.\nabla p_{\vartheta}(\mathsf{x}_{n}|\mathsf{Y}_{1:n-1})\right\vert _{\vartheta=\vartheta_{n}}$
respectively. It has been shown in \cite{LeGlandMevel97} that the
solution of RML converges to the true ML estimator without any loss
of efficiency. For more details on the convergence of RML for HMMs
we refer the reader to \cite{LeGlandMevel97}.



\subsection{On-line Expectation-Maximization (EM) \label{subsec:em}}

We begin this section with a brief description of Expectation-Maximization
(EM) \cite{Demster77} and then present its on-line implementation.
EM is an iterative off-line algorithm for learning $\vartheta_{\ast}$,
which consists of repeating a two step procedure given a batch of
$T$\ observations. Let $p$ be the (off-line) iteration index. The
first step, the expectation or E-step, computes 
\begin{equation}
Q(\vartheta_{p},\vartheta)=\int\log p_{\vartheta}(\mathsf{x}_{1:T},\mathsf{Y}_{1:T})p_{\vartheta_{p}}(\mathsf{x}_{1:T}|\mathsf{Y}_{1:T})d\mathsf{x}_{1:T}.\label{eq:Qfunction}
\end{equation}
 The second step is the maximization or M-step that updates the parameter
$\vartheta_{p}$, 
\begin{equation}
\vartheta_{p+1}=\arg\max\quad Q(\vartheta_{p},\vartheta)\label{eq:maximizationEM}
\end{equation}
 Upon the completion of an E and M step, the likelihood surface is
ascended, i.e. $p_{\vartheta_{p+1}}(\mathsf{Y}_{1:T})\geq p_{\vartheta_{p}}(\mathsf{Y}_{1:T})$
\cite{Demster77}. When $p_{\theta}(\mathsf{x}_{1:T},\mathsf{Y}_{1:T})$
is in the exponential family, which is the case of linear Gaussian
state-space models, this procedure can be implemented exactly. Then
the E-step is equivalent to computing a summary statistic of the form
\begin{equation}
\mathcal{S}_{T}^{\vartheta_{p}}=\frac{1}{T}\int\left(\sum_{n=1}^{T}s_{n}\left(\mathsf{x}{}_{n-1:n},\mathsf{Y}_{n}\right)\right)p_{\vartheta_{p}}(\mathsf{x}_{1:T}|\mathsf{Y}_{1:T})d\mathsf{x}_{1:T}.\label{eq:summaryEM}
\end{equation}
 where $s_{n}:\mathbb{R}^{d_{x}}\times\mathbb{R}^{d_{x}}\times\mathbb{R}^{d_{y}}\rightarrow\mathbb{R}^{\kappa}$.
In addition, the maximizing argument of $Q(\vartheta_{p},\vartheta)$
can be characterized in this case explicitly through a suitable function
$\Lambda:\mathbb{R}^{\kappa}\rightarrow\Theta$, i.e. 
\begin{equation}
\vartheta_{p}=\Lambda\left(\mathcal{S}_{T}^{\vartheta_{p}}\right).\label{eq:maximiEM}
\end{equation}
 Note that in the usual EM setup one has to compute (\ref{eq:summaryEM})
for every iteration $p$ of the algorithm.

It is also possible to propose an on-line version of the EM algorithm.
This was originally proposed for finite state-space and linear Gaussian
models in \cite{elliott-online-em,elliott-ijacsp,ford-thesis} and
for exponential family models in \cite{cappe-on-line-em,cappe-on-line-em-hmm}.
In the online implementation of the EM, running averages of the sufficient
statistics are computed \cite{ford-thesis,elliott-online-em,cappe-on-line-em}.
Let $\{\vartheta_{m}\}_{1\leq m\leq n}$ be the sequence of parameter
estimates of the online EM algorithm computed sequentially based on
$\mathsf{Y}_{1:n-1}$. When $\mathsf{Y}_{n}$ is received, we compute
\begin{equation}
\begin{tabular}{l}
 \ensuremath{\mathcal{S}_{n}=\gamma_{n}\text{ }\int s_{n}\left(\mathsf{x}{}_{n-1:n}\right)p_{\vartheta_{1:n}}(\mathsf{x}{}_{n-1:n}|\mathsf{Y}_{1:n})d\mathsf{x}{}_{n-1:n}}\\
\ensuremath{+\left(1-\gamma_{n}\right)\sum_{m=1}^{n-1}(\prod\limits _{i=m+1}^{n-1}\left(1-\gamma_{i}\right))\gamma_{m}\int}\ensuremath{s_{m}\left(\mathsf{x}_{m-1:m}\right)p_{\vartheta_{1:n}}(\mathsf{x}{}_{m-1:m}|\mathsf{Y}_{1:n})d\mathsf{x}{}_{m-1:m},}
\end{tabular}\label{eq:suffStatOnline}
\end{equation}
 where the subscript $\vartheta_{1:n}$ on $p_{\vartheta_{1:n}}(\mathsf{x}_{1:T}|\mathsf{Y}_{1:n})$
indicates that the posterior density is being computed sequentially
using the parameter $\vartheta_{m}$ at time $m\leq n$. The step
sizes $\left\{ \gamma_{n}\right\} _{n\geq1}$ need to satisfy $\sum_{n}\gamma_{n}=\infty$
and $\sum_{n}\gamma_{n}^{2}<\infty$ as in the RML case. For the M-step
one uses the same maximization step (\ref{eq:maximiEM}) used in the
batch version 
\begin{equation}
\vartheta_{n+1}=\Lambda\left(\mathcal{S}_{n}\right).
\end{equation}
 The recursive calculation of $\mathcal{S}_{n}$ can be achieved by
setting $V_{1}\left(\mathsf{x}_{0}\right)=0$ and computing 
\begin{eqnarray}
V_{n}\left(\mathsf{x}_{n}\right) & = & \int\left\{ \gamma_{n}\text{ }s_{n}\left(\mathsf{x}_{n-1},\mathsf{x}_{n}\right)+\left(1-\gamma_{n}\right)\text{ }V_{n-1}\left(\mathsf{x}_{n-1}\right)\right\} \nonumber \\
 & \times & p_{\vartheta_{1:n}}\left(\left.\mathsf{x}_{n-1}\right\vert \mathsf{Y}_{1:n-1},\mathsf{x}_{n}\right)d\mathsf{x}_{n-1}\label{eq:updateStatonline}
\end{eqnarray}
 and 
\begin{equation}
\mathcal{S}_{n}=\int V_{n}\left(\mathsf{x}_{n}\right)p_{\vartheta_{1:n}}(\mathsf{x}_{n}|\mathsf{Y}_{1:n})d\mathsf{x}_{n}.\label{eq:updateStatonline2}
\end{equation}
 For finite state-space and linear Gaussian models, all the quantities
appearing in this algorithm can be calculated exactly \cite{ford-thesis,elliott-online-em,cappe-on-line-em-hmm}.

\section{Distributed RML derivation}


Let $\theta_{n}=\{\theta_{n}^{i,j}\}_{(i,j)\in\mathcal{E}}$ be the
estimate of the true parameter $\theta_{\ast}$\ given the available
data $Y_{1:n-1}$. Consider an arbitrary node $r$ and assume it controls
edge $(r,j)$. At time $n$, we assume the following quantities are
available: $(\dot{\mu}_{n-1}^{r,j}=\left.\nabla_{\theta^{r,j}}\mu_{n-1}^{r}\right\vert _{\theta=\theta_{n}},\left.\mu_{n-1}^{r}\right\vert _{\theta=\theta_{n}},\Sigma_{n-1}^{r})$.
The first of these quantities is the derivative of the conditional
mean of the hidden state at node $r$ given $Y_{1:n-1}$, i.e. $\left.\nabla_{\theta^{r,j}}\int x_{n-1}^{r}p_{\theta}^{r}(x_{n-1}^{r}|Y_{1:n-1})dx_{n-1}^{r}\right\vert _{\theta=\theta_{n}}$.
This quantity is a function of the localization parameter $\theta_{n}$.
$\Sigma_{n-1}^{r}$\ is the variance of the distribution $\left.p_{\theta}^{r}(x_{n-1}^{r}|Y_{1:n-1})\right\vert _{\theta=\theta_{n}}$\ and
is independent of the localization parameter. The log-likelihood in
(\ref{eq:almostOnlineDistributedRML}) evaluates to: 
\begin{align*}
 & \log p_{\theta}^{r}(Y_{n}|Y_{1:n-1})=-\frac{1}{2}\underset{i\in\mathcal{V}}{\sum}(Y_{n}^{i}-C_{n}^{i}\theta^{r,i})^{T}R_{n}^{i}{}^{-1}(Y_{n}^{i}-C_{n}^{i}\theta^{r,i})\\
 & -\frac{1}{2}\mu_{n|n-1}^{r}{}^{\text{T}}(\Sigma_{n|n-1}^{r})^{-1}\mu_{n|n-1}^{r}+\frac{1}{2}(z_{n}^{r})^{\text{T}}(M_{n}^{r})^{-1}z_{n}^{r}+const
\end{align*}
 where all $\theta$ independent terms have been lumped together in
the term `const'.\ (Refer to Algorithm \ref{alg:distFilter}\ for
the definition of the quantities in this expression.) Differentiating
this expression w.r.t. $\theta^{r,j}$\ yields 
\begin{align*}
\nabla_{\theta^{r,j}}\log p_{\theta}^{r}(Y_{n}|Y_{1:n-1}) & =-(\nabla_{\theta^{r,j}}\mu_{n|n-1}^{r})^{\text{T}}(\Sigma_{n|n-1}^{r})^{-1}\mu_{n|n-1}^{r}\\
 & +(\nabla_{\theta^{r,j}}z_{n}^{r})^{\text{T}}(M_{n}^{r})^{-1}z_{n}^{r}+\underset{i\in\mathcal{V}}{\sum}(\nabla_{\theta^{r,j}}\theta^{r,i})^{\text{T}}(C_{n}^{i})^{\text{T}}(R_{n}^{i})^{-1}(Y_{n}^{i}-C_{n}^{i}\theta^{r,i}).
\end{align*}
 (\ref{eq:almostOnlineDistributedRML}) requires $\nabla_{\theta^{r,j}}\log p_{\theta}^{r}(Y_{n}|Y_{1:n-1})$
to be evaluated at $\theta=\theta_{n}$. Using the equations (\ref{eq:distrKalmanSecond})-(\ref{eq:distrKalmanLast})
and the assumed knowledge of $(\dot{\mu}_{n-1}^{r,j},\left.\mu_{n-1}^{r}\right\vert _{\theta=\theta_{n}},\Sigma_{n-1}^{r})$\ we
can evaluate the derivatives on the right-hand side of this expression:
\begin{align}
\dot{\mu}_{n|n-1}^{r,j} & =\left.\nabla_{\theta^{r,j}}\mu_{n|n-1}^{r}\right\vert _{\theta=\theta_{n}}=A_{n}\dot{\mu}_{n-1}^{r,j},\label{eq:gradientsPropagate1}\\
\dot{z}_{n}^{r,j} & =\left.\nabla_{\theta^{r,j}}z_{n}^{r}\right\vert _{\theta=\theta_{n}}=(\Sigma_{n|n-1}^{r})^{-1}\dot{\mu}_{n|n-1}^{r,j}-\underset{i\in\mathcal{V}}{\sum}(C_{n}^{i})^{\text{T}}(R_{n}^{i})^{-1}C_{n}^{i}\left.\nabla_{\theta^{r,j}}\theta^{r,i}\right\vert _{\theta=\theta_{n}},\label{eq:gradientsPropagate2}\\
\dot{\mu}_{n}^{r,j} & =\left.\nabla_{\theta^{r,j}}\mu_{n}^{r}\right\vert _{\theta=\theta_{n}}=(M_{n}^{r})^{-1}\dot{z}_{n}^{r,j}.\label{eq:gradientsPropagate3}
\end{align}
 Using property (\ref{eq:sumThetas}) we note that for the set of
vertices $i$ for which the path from $r$ to $i$ includes edge $(r,j)$,
$\nabla_{\theta^{r,j}}\theta^{r,i}=I$ (the identity matrix) whereas
for the rest $\nabla_{\theta^{r,j}}\theta^{r,i}=0$. For all the nodes
$i$ for which $\nabla_{\theta^{r,j}}\theta^{r,i}=I$, let them form
a sub tree $(\mathcal{V}_{rj}^{\prime},\mathcal{E}_{rj}^{\prime})$
branching out from node $j$ away from node $r$. Then the last sum
in the expression for $\left.\nabla_{\theta^{r,j}}\log p_{\theta}^{r}(Y_{n}|Y_{1:n-1})\right\vert _{\theta=\theta_{n}}$
evaluates to, 
\[
\underset{i\in\mathcal{V}_{rj}^{\prime}}{\sum}(C_{n}^{i})^{\text{T}}(R_{n}^{i})^{-1}(Y_{n}^{i}-C_{n}^{i}\theta_{n}^{r,i})=\dot{m}_{n,K}^{j,r}-\ddot{m}_{n,K}^{j,r},
\]
 where messages $(\dot{m}_{n,K}^{j,r},\ddot{m}_{n,K}^{j,r})$ were
defined in Algorithms \ref{alg:distFilter}. Similarly, we can write
the sum in the expression for $\dot{z}_{n}^{r,j}$ as $m_{n,K}^{j,r}$
(again refer to Algorithms \ref{alg:distFilter}) to obtain 
\begin{equation}
\dot{z}_{n}^{r,j}=(\Sigma_{n|n-1}^{r})^{-1}\dot{\mu}_{n|n-1}^{r,j}-m_{n,K}^{j,r}.\label{eq:gradientsPropagate22}
\end{equation}
 To conclude, the approximations to $(\dot{\mu}_{n}^{r,j}=\left.\nabla_{\theta^{r,j}}\mu_{n}^{r}\right\vert _{\theta=\theta_{n+1}},\left.\mu_{n}^{r}\right\vert _{\theta=\theta_{n+1}},\Sigma_{n}^{r})$
for the subsequent RML iteration, i.e. (\ref{eq:almostOnlineDistributedRML})
at time $n=n+1,$ are given by 
\[
\dot{\mu}_{n}^{r,j}=(M_{n}^{r})^{-1}\dot{z}_{n}^{r,j}
\]
 while $(\left.\mu_{n}^{r}\right\vert _{\theta=\theta_{n+1}},\Sigma_{n}^{r})$
are given by (\ref{eq:distrKalmanSecond})-(\ref{eq:distrKalmanLast}).
The approximation to $\left.\nabla_{\theta^{r,j}}\mu_{n}^{r}\right\vert _{\theta=\theta_{n+1}}$
follows from differentiating (\ref{eq:distrKalmanLast}). $(\dot{\mu}_{n}^{r,j},\left.\mu_{n}^{r}\right\vert _{\theta=\theta_{n+1}})$
are only approximations because they are computed using the previous
values of the parameters, i.e. $\theta_{1:n}$.


\section{Distributed EM derivation}

For the off-line EM approach, once a batch of $T$ observations have
been obtained, each node $r$ of the network that controls an edge
will execute the following E and M step iteration $n$, 
\begin{align*}
Q^{r}(\theta_{p},\theta) & =\int\log p_{\theta}^{r}(x_{1:T}^{r},Y_{1:T})p_{\theta_{p}}^{r}(x_{1:T}^{r}|Y_{1:T})dx_{1:T}^{r},\\
\theta_{p+1}^{r,j} & =\arg\underset{\theta^{r,j}\in\Theta}{\max}Q^{r}(\theta_{p},(\theta^{r,j},\{\theta^{e},e\in\mathcal{E}\backslash(r,j)\})),
\end{align*}
 where it is assumed that node $r$ controls edge $(r,j)$. The quantity
$p_{\theta_{p}}^{r}(x_{1:T}^{r}|Y_{1:T})$\ is the joint distribution
of the hidden states at node $r$ given all the observations of the
network from time $1$ to $T$ and is given up to a proportionality
constant, 
\[
p_{\theta_{p}}^{r}(x_{1:T}^{r})p_{\theta_{p}}^{r}(Y_{1:T}|x_{1:T}^{r})=\prod\limits _{n=1}^{T}f_{n}^{}(x_{n}^{r}|x_{n-1}^{r})p_{\theta_{p}}^{r}(Y_{n}|x_{n}^{r}),
\]
 where $p_{\theta_{p}}^{r}(Y_{n}|x_{n}^{r})$\ was defined in (\ref{eq:collabLikelihood}).
Note that $p_{\theta_{p}}^{r}(x_{1:T}^{r},Y_{1:T})$\ (and hence
$p_{\theta_{p}}^{r}(x_{1:T}^{r}|Y_{1:T})$) is a function of $\theta_{p}=\{\theta_{p}^{i,i^{\prime}}\}_{(i,i^{\prime})\in\mathcal{E}}$\ and
not just $\theta_{p}^{r,j}$. Also, the $\theta$-dependence of $p_{\theta}^{r}(x_{1:T}^{r},Y_{1:T})$\ arises
through the likelihood term only as $p_{\theta}^{r}(x_{1:T}^{r})$\ is
$\theta$ independent. Note that 
\begin{align*}
\sum_{v\in\mathcal{V}}\log g_{n}^{v}(Y_{n}^{v}|x_{n}^{r}+\theta^{r,v}) & =\sum_{v\in\mathcal{V}}c_{n}^{v}-\frac{1}{2}\sum_{v\in\mathcal{V}}(Y_{n}^{v}-C_{n}^{v}\theta^{r,v})^{\text{T}}(R_{n}^{v})^{-1}(Y_{n}^{v}-C_{n}^{v}\theta^{r,v})\\
 & +(x_{n}^{r})^{\text{T}}\sum_{v\in\mathcal{V}}(C_{n}^{v})^{\text{T}}(R_{n}^{v})^{-1}(Y_{n}^{v}-C_{n}^{v}\theta^{r,v})-\frac{1}{2}(x_{n}^{r})^{\text{T}}\left[\sum_{v\in\mathcal{V}}(C_{n}^{v})^{\text{T}}(R_{n}^{v})^{-1}C_{n}^{v}\right]x_{n}^{r}
\end{align*}
 where $c_{n}^{v}$ is a constant independent of $\theta$. Taking
the expectation w.r.t. $p_{\theta_{n}}^{r}(x_{n}^{r}|Y_{1:T})$ gives
\begin{align*}
\int\log p_{\theta}^{r}(Y_{n}|x_{n}^{r})p_{\theta_{p}}^{r}(x_{n}^{r}|Y_{1:T})dx_{n}^{r} & =-\frac{1}{2}\sum_{v\in\mathcal{V}}\left[(Y_{n}^{v}-C_{n}^{v}\theta^{r,v})^{\text{T}}(R_{n}^{v})^{-1}(Y_{n}^{v}-C_{n}^{v}\theta^{r,v})\right]\\
 & -(\mu_{n|T}^{r})^{\text{T}}\sum_{v\in\mathcal{V}}(C_{n}^{v})^{\text{T}}(R_{n}^{v})^{-1}C_{n}^{v}\theta^{r,v}+\text{const}
\end{align*}
 where all terms independent of $\theta^{r,j}$\ have been lumped
together as 'const' and $\mu_{n|T}^{r}$\ is the mean of $x_{n}^{r}$\ under
$p_{\theta_{p}}^{r}(x_{n}^{r}|Y_{1:T})$. Taking the gradient w.r.t.
$\theta^{r,j}$ and following the steps in the derivation of the distributed
RML we obtain 
\[
\nabla_{\theta^{r,j}}\int\log p_{\theta}^{r}(Y_{n}|x_{n}^{r})p_{\theta_{p}}^{r}(x_{n}^{r}|Y_{1:T})dx_{n}^{r}=\dot{m}_{n,K}^{j,r}-\ddot{m}_{n,K}^{j,r}-(m_{n,K}^{j,r})^{\text{T}}\mu_{n|T}^{r}
\]
 where $(m_{n,K}^{j,r},\dot{m}_{n,K}^{j,r},\ddot{m}_{n,K}^{j,r})$
is defined in (\ref{eq:message1})-(\ref{eq:message3}). Only $\ddot{m}_{n,K}^{j,r}$\ is
a function of $\theta^{r,j}$. Now to perform the M-step, we solve
\[
\left(\sum_{n=1}^{T}m_{n,K}^{j,r}\right)\theta^{r,j}=\sum_{n=1}^{T}\left(\dot{m}_{n,K}^{j,r}-(m_{n,K}^{j,r})^{\text{T}}\mu_{n|T}^{r}-\sum\limits _{j'\in\text{ne}(j)\setminus\{r\}}\ddot{m}_{n,K-1}^{j',j}\right).
\]
 Note that $\theta^{r,j}$\ can be recovered by standard linear algebra
and so far $\theta^{r,j}$\ is solved by quantities available locally
to node $r$ and $j$. One can use the fact that $\sum\limits _{j'\in\text{ne}(j)\setminus\{r\}}\ddot{m}_{n,K-1}^{j',j}=\ddot{m}_{n,K}^{j,r}-\ddot{m}_{n,1}^{j,r}$
to so that the M-step can be performed with quantities available locally
to node $r$ only. Recall that $\sum_{n=1}^{T}\mu_{n|T}^{r}=\int\left(\sum_{n=1}^{T}x_{n}^{r}\right)p_{\theta_{p}}^{r}(x_{1:T}^{r}|Y_{1:T})dx_{1:T}^{r}.$
This implies directly that three summary statistics are needed for
node $r$ to update $\theta^{r,j}$. These should be defined using:
\[
s_{n,1}^{r,j}(x_{n}^{r},Y_{n})=(m_{n,K}^{j,r})^{\text{T}}x_{n}^{r},\: s_{n,2}^{r,j}(x_{n}^{r},Y_{n})=m_{n,K}^{j,r},\: s_{n,3}^{r,j}(x_{n}^{r},Y_{n})=\dot{m}_{n,K}^{j,r}-\ddot{m}_{n,K}^{j,r}+\ddot{m}_{n,1}^{j,r},
\]
 where $s_{n,1}^{r}$, $s_{n,3}^{r}$ are each functions of $x_{n}^{r}$
and $Y_{n}$ via $\mu_{n|T}^{r}$ and $\dot{m}_{n,K}^{j,r}-\ddot{m}_{n,K}^{j,r}+\ddot{m}_{n,1}^{j,r}$
respectively. The summary statistics can be written in the form of
(\ref{eq:summaryEM}) as follows: 
\begin{eqnarray*}
\mathcal{S}_{T,1}^{r,j^{,\theta_{p}}} & = & \frac{1}{T}\int\left(\sum_{n=1}^{T}(m_{n,K}^{j,r})^{\text{T}}x_{n}^{r}\right)p_{\theta_{p}}^{r}(x_{1:T}^{r}|Y_{1:T})dx_{1:T}^{r},\\
\mathcal{S}_{T,2}^{r,j^{\theta_{p}}} & = & \frac{1}{T}\sum_{n=1}^{T}m_{n,K}^{j,r},\quad\mathcal{S}_{T,3}^{r,j^{\theta_{p}}}=\frac{1}{T}\sum_{n=1}^{T}\left(\dot{m}_{n,K}^{j,r}-\ddot{m}_{n,K}^{j,r}+\ddot{m}_{n,1}^{j,r}\right),
\end{eqnarray*}
 and the M-step function becomes $\Lambda(s_{1},s_{2},s_{3})=s_{2}^{-1}\left(s_{3}-s_{1}\right),$
where $s_{1},s_{2},s_{3}$ correspond to each of the three summary
statistics. Note that $\Lambda$ is the same function for every node.

We will now proceed to the on-line implementation. Let at time $n$
the estimate of the localization parameter be $\theta_{n}$. Following
the description of Section \ref{subsec:em}, for every $r\in\mathcal{V}$
and $(r,j)\in\mathcal{E}$, let $\mathcal{S}_{n,1}^{r,j}$, $\mathcal{S}_{n,2}^{r,j}$,
$\mathcal{S}_{n,3}^{r,j}$ be the running averages (w.r.t $n$) for
$\mathcal{S}_{T,1}^{r,j^{\theta_{p}}}$, $\mathcal{S}_{T,2}^{r,j^{\theta_{p}}}$
and $\mathcal{S}_{T,3}^{r,j^{\theta_{p}}}$ respectively. The recursions
for $\mathcal{S}_{n,2}^{r,j}$, $\mathcal{S}_{n,3}^{r,j}$ are trivial:
\[
\mathcal{S}_{n,2}^{r,j}=\gamma_{n}^{r}m_{n,K}^{j,r}+\left(1-\gamma_{n}^{r}\right)\mathcal{S}_{n-1,2}^{r,j},\:\mathcal{S}_{n,3}^{r,j}=\gamma_{n}^{r}(\dot{m}_{n,K}^{j,r}-\ddot{m}_{n,K}^{j,r}+\ddot{m}_{n,1}^{j,r})+\left(1-\gamma_{n}^{r}\right)\mathcal{S}_{n-1,3}^{r,j},
\]
 where $\left\{ \gamma_{n}^{r}\right\} _{n\geq1}$ needs to satisfy
$\sum_{n\geq\text{1}}\gamma_{n}^{r}=\infty$ and $\sum_{n\geq1}\left(\gamma_{n}^{r}\right)^{2}<\infty$.
For $\mathcal{S}_{n,1}^{r,j}$, we will use (\ref{eq:updateStatonline})-(\ref{eq:updateStatonline2}).
We first set $V_{0}^{r,j}\left(x_{0}^{r}\right)=0$ and define the
recursion 
\begin{equation}
V_{n}^{r,j}\left(x_{n}^{r}\right)=\gamma_{n}^{r}(m_{n,K}^{j,r})^{\text{T}}x_{n}^{r}+\left(1-\gamma_{n}^{r}\right)\int\text{ }V_{n-1}^{r,j}\left(x_{n-1}^{r}\right)p_{\theta_{1:n}}^{r}\left(\left.x_{n-1}^{r}\right\vert Y_{1:n-1},x_{n}^{r}\right)dx_{n-1}^{r}.\label{eq:VrecursionApp}
\end{equation}
 Using standard manipulations with Gaussians we can derive that $p_{\theta_{1:n}}^{r}\left(\left.x_{n-1}^{r}\right\vert Y_{1:n-1},x_{n}^{r}\right)$
is itself a Gaussian density with mean and variance denoted by $\tilde{\mu}_{n}^{r}(x_{n}),\widetilde{\Sigma}_{n}^{r}$
respectively, where 
\[
\tilde{\Sigma}_{n}^{r}=\left(\Sigma_{n-1}^{r}+A_{n}^{\text{T}}Q_{n}^{-1}A_{n}\right)^{-1},\:\widetilde{\mu}_{n}^{r}(x_{n})=\tilde{\Sigma}_{n}^{r}\left(\left(\Sigma_{n-1}^{r}\right)^{-1}\mu_{n-1}^{r}+A_{n}^{\text{T}}Q_{n}^{-1}x_{n}\right).
\]
 It is then evident that (\ref{eq:VrecursionApp}) becomes $V_{n}^{r,j}\left(x_{n}^{r}\right)=H_{n}^{r,j}x_{n}^{r}+h_{n}^{r,j}$,
with: 
\begin{eqnarray*}
H_{n}^{r,j} & = & \gamma_{n}^{r}(m_{n,K}^{j,r})^{\text{T}}+(1-\gamma_{n}^{r})H_{n-1}^{r,j}\left(\tilde{\Sigma}_{n}^{r}\right)^{-1}A_{n}^{\text{T}}Q_{n}^{-1},\\
h_{n}^{r,j} & = & (1-\gamma_{n}^{r})\left(H_{n-1}^{r,j}\left(\tilde{\Sigma}_{n}^{r}\right)^{-1}\left(\Sigma_{n-1}^{r}\right)^{-1}\mu_{n-1}^{r}+h_{n-1}^{r,j}\right),
\end{eqnarray*}
 where $H_{0}^{r,j}=0$ and $h_{0}^{r,j}=0$. Finally, the recursive
calculation of $\mathcal{S}_{n,1}^{r,j}$ is achieved by computing
\[
\mathcal{S}_{n,1}^{r,j}=\int V_{n}^{r,j}\left(x_{n}^{r}\right)p_{\theta_{0:n}}^{r}(x_{n}^{r}|Y_{0:n})dx_{n}^{r}=H_{n}^{r,j}\mu_{n}^{r}+h_{n}^{r,j}.
\]
 Again all the steps are performed locally at node $r$, which can
update parameter $\theta^{r,j}$ using $\theta_{n+1}^{r,j}=\Lambda(\mathcal{S}_{n,1}^{r,j},\mathcal{S}_{n,2}^{r,j},\mathcal{S}_{n,3}^{r,j}).$

\section*{Acknowledgement}

N. Kantas was supported by the Engineering and Physical Sciences Research
Council programme grant on Control For Energy and Sustainability (EP/G066477/1).
S.S. Singh's research is partly funded by the Engineering and Physical
Sciences Research Council under the First Grant Scheme (EP/G037590/1).

\bibliographystyle{IEEEtran}

\addcontentsline{toc}{section}{\refname}\bibliographystyle{IEEEtran}
\bibliography{First}

\end{document}